\documentclass[review]{elsarticle}

 \usepackage[a4paper, total={5.5in,8in}]{geometry}
\usepackage[utf8]{inputenc}
 \usepackage{lipsum}
 \usepackage{amsmath,amssymb}
\usepackage{mathtools}
\usepackage{color}
\usepackage{xcolor, hyperref}
\hypersetup{
    colorlinks = true,
    linkcolor = red,
    linkbordercolor = {white}
}

\usepackage{graphicx}
\usepackage{amsthm}
\usepackage{ulem}
\usepackage{caption}
\usepackage{subcaption}

\newcommand{\usigma}{\sigma}
\newcommand{\uvarepsilon}{\varepsilon}
\newcommand{\victor}[1]{\textcolor{black}{#1}}
\newcommand{\elisabeth}[1]{\textcolor{black}{#1}}
\newcommand{\new}[1]{\textcolor{black}{#1}}

\DeclareMathOperator{\diag}{diag}
\DeclareMathOperator{\Tr}{Tr}

\begin{document}

\begin{frontmatter}
\title{An unfitted radial basis function generated finite difference method applied to thoracic diaphragm simulations}

\author[igor]{Igor Tominec}
\ead[igor]{igor.tominec@it.uu.se}
\address[igor]{Uppsala University, Department of Information Technology, Division of Scientific Computing, Sweden}

\author[pierre]{Pierre-Frédéric Villard}
\ead[pierre]{pierrefrederic.villard@loria.fr}
\address[pierre]{Université de Lorraine, CNRS, Inria, LORIA, Nancy, France}

\author[igor]{Elisabeth Larsson}
\ead[elisabeth]{elisabeth.larsson@it.uu.se}

\author[victor]{Víctor Bayona}
\ead[victor]{vbayona@math.uc3m.es}
\address[victor]{Departamento de Matemáticas, UC3M, Leganés 28911, Spain}

\author[nicola,nicola2]{Nicola Cacciani}
\ead[nicola]{nicola.cacciani@ki.se}
\address[nicola]{Department of Physiology and Pharmacology, Karolinska Institutet,Stockholm, Sweden}
\address[nicola2]{Department of Clinical Neuroscience, Clinical Neurophysiology, Karolinska Institutet, Stockholm, Sweden}

\begin{abstract}
    The thoracic diaphragm is the muscle that drives the respiratory cycle of a human being. \elisabeth{Using a system of partial differential 
    equations (PDEs) that models linear elasticity} we compute displacements and stresses 
    in a two-dimensional cross section of the diaphragm in its contracted state. 
    The boundary data consists of a mix of displacement and traction conditions. 
    If these are imposed as they are, and the conditions are not compatible, this leads to reduced smoothness of the solution. 
    Therefore, the boundary data is first smoothed using the least-squares radial basis function generated finite difference (RBF-FD) framework.
    Then the boundary conditions are reformulated as a Robin boundary condition with smooth coefficients. 
    The same framework is also used to approximate the boundary curve of the diaphragm cross section based on data obtained from a slice of a 
    computed tomography (CT) scan. To solve the 
    PDE we employ the unfitted least-squares RBF-FD method. This makes it easier to handle the geometry of the diaphragm, 
    which is thin and non-convex.
    We show numerically that our solution converges with high-order 
    towards a finite element solution evaluated on a fine grid.
    Through this simplified numerical model we also gain an insight into the challenges 
    associated with the diaphragm geometry and the boundary conditions before approaching a more complex three-dimensional model. 
\end{abstract}
\end{frontmatter}


\section{Introduction}
During the 2020 covid-19 pandemic we have all been made aware that intensive care units (ICU) have limited capacity with respect to the number of patients that can be cared for at the same time. The WHO report on covid-19 in China~\cite{WHO} indicates that covid-19 patients with severe symptoms need 3--6 weeks in ICU.
Patients with severe respiratory symptoms are put under mechanical ventilation to save their lives. At the same time, the mechanical ventilation has adverse effects, such as ventilator induced diaphragmatic dysfunction (VIDD)~\cite{Cacciani}, which prolongs the ICU time.  Hence, improving mechanical ventilation can have a significant impact on ICU patient turnover. 

This work is part of the INVIVE project\footnote{\url{https://www.it.uu.se/research/scientific_computing/project/rbf/biomech}}, 
where we aim to create a mechanically ventilated virtual patient~\cite{KZECH15}
on whom we can perform tests with different ventilation strategies and counter measures against VIDD~\cite{Jaber17}. With the virtual patient, we can also vary the gender, physiology, age and potential injuries that affect the response. Using simulations allows to create a controlled computer-based environment in a way that is not possible in a clinical setting.

The diaphragm is the main respiratory muscle, and the present focus of our study. The diaphragm is located between the thoracic and abdominal cavities. It has two domes, and is attached to the lower ribs, the spine, and the sternum. 
During mechanical ventilation, the normal action of the diaphragm, where inhalation follows the contraction of the muscle, is reversed, i.e., 
the muscles become passive and the air is pumped into the lungs by the ventilator. As air is entering the lungs with a positive pressure, the muscle fibres 
are instead extended. This sudden and extreme mechanical perturbation is the trigger of a chain of biological events causing 
the progression of VIDD.


Numerical simulation of the biomechanical action of the diaphragm during respiration or ventilation is a challenging problem. The shape of the diaphragm is non-trivial and there are gradual transitions between muscle and tendon with very different material response. Imperfect data can be extracted from medical images, and can then be converted into a geometry representation~\cite{invive}. The specific challenges of constructing a smooth geometry representation are addressed in a 
forthcoming paper~\cite{Larsson20}. 



In this paper, we model the diaphragm using a linear elastic PDE. This is a simplification and we plan to develop 
a more realistic tissue model as a part of our future work.

The boundary conditions for the elastic PDE system are given by a combination of traction boundary conditions, 
which contains first derivatives of the displacement, and (time-dependent) Dirichlet boundary conditions for the displacement, 
where the diaphragm is attached. These are mixed boundary conditions, in general not fully compatible, 
leading to reduced regularity of the solution even when we expect a smooth solution 
from a physiological perspecitve. 
Therefore, we smooth the boundary data as well as the transition between the traction boundary data and Dirichlet 
boundary data before using it in the PDE solver. 


We use the unfitted radial basis function generated finite difference method in the least-squares setting 
(unfitted RBF-FD-LS)~\cite{tominec2020unfitted} to solve the diaphragm problem. A benefit of using the unfitted RBF-FD-LS method is that the PDE problem 
is solved on an extended domain (see Figure \ref{fig:method:diaphragm_interpNodes}). This simplifies node generation 
as the node placement is then independent of the geometry.
Another benefit of the unfitted setting is a smaller approximation error near the boundaries \cite{tominec2020unfitted}, where the stencils 
are typically highly skewed when using conventional RBF-FD methods such as the fitted RBF-FD-LS method \cite{ToLaHe20} 
and the collocation RBF-FD method \cite{FoFly15book}. 

As model geometry we use a two-dimensional cross section of the diaphragm.
To investigate the properties of the problem we use a combination of data from medical images and knowledge about the physiology of the diaphragm 
expressed in terms of boundary conditions. In particular, we construct two benchmark problems: (i) a pure Dirichlet case, (ii) 
a case with mixed boundary conditions (Dirichlet + traction), rephrased as a Robin condition with smooth coefficients.
With the experiments performed for these benchmarks, we aim to answer the following questions:
\begin{itemize}
\item Can we solve the benchmarks problems with unfitted RBF-FD-LS? Are there specific numerical challenges to be noted?
\item Can we achieve high-order convergence to the solution of the elastic PDE when the imposed boundary data and the geometry is smooth?
\item How is the performance and accuracy of the RBF-FD-LS solver affected by the type of boundary conditions that are imposed?
\item How does the RBF-FD-LS solver compare with a basic FEM solver? Do the solvers give similar solutions?
\end{itemize}

\new{Other authors have also modeled the diaphragm numerically.} The most advanced diaphragm simulations in the literature can be found in a series of publications by Ladjal et al.~\cite{Ladjal13,Ladjal15,Ladjal15b,Giroux17,Ladjal19}. The application focus is to track the motion of lung tumours during respiration for radiotherapy purposes. The finite element models that are employed are highly elaborate, taking into account features such as patient-specific lung compliance in order to define the constitutive law. 
The simulation times reported are far from real-time capability. In these studies, the diaphragm is not the main target and it  is not well-resolved in the thickness direction due to the high aspect ratio.
%
Other relevant, but less detailed, diaphragm simulations can be found in~\cite{Pato11}, where the diaphragm is discretized using shell elements to compare healthy and pathological situations, and in~\cite{Fuerst14}, where the whole region under the lungs, including the diaphragm, forms one region in the simulation.

The outline of the paper is as follows: Section~\ref{section:diabehavior}
describes the expected behavior of the diaphragm during respiration and the relation of the two-dimensional 
geometry to the real three-dimensional diaphragm geometry. 
Then the linear elasticity model problem is defined in Section~\ref{sec:linearelasticity}. 
An RBF-FD algorithm for computing differentiation and evaluation matrices is described in Section~\ref{sec:rbffd}. 
In Section~\ref{sec:unfitted_discretization} these matrices are then used in the least-squares unfitted RBF-FD setting 
to discretize the linear elasticity equations.
Section ~\ref{sec:smooth} discusses how the boundary of the two-dimensional diaphragm cross-section, and the boundary conditions, are smoothed. 
In Section ~\ref{sec:experiments:dirichlet} and Section~\ref{sec:experiments:Robin} we compute the solution to the linear elasticity equations 
for the two benchmark problems, and evaluate the convergence numerically. 
Finally, Section~\ref{sec:conclusions} contains the conclusions.

\section{The expected behavior of the thoracic diaphragm}
\label{section:diabehavior}
\begin{figure}
  \centering
  \begin{subfigure}[b]{0.48\textwidth}
    \centering
    \includegraphics[height=0.55\textwidth]{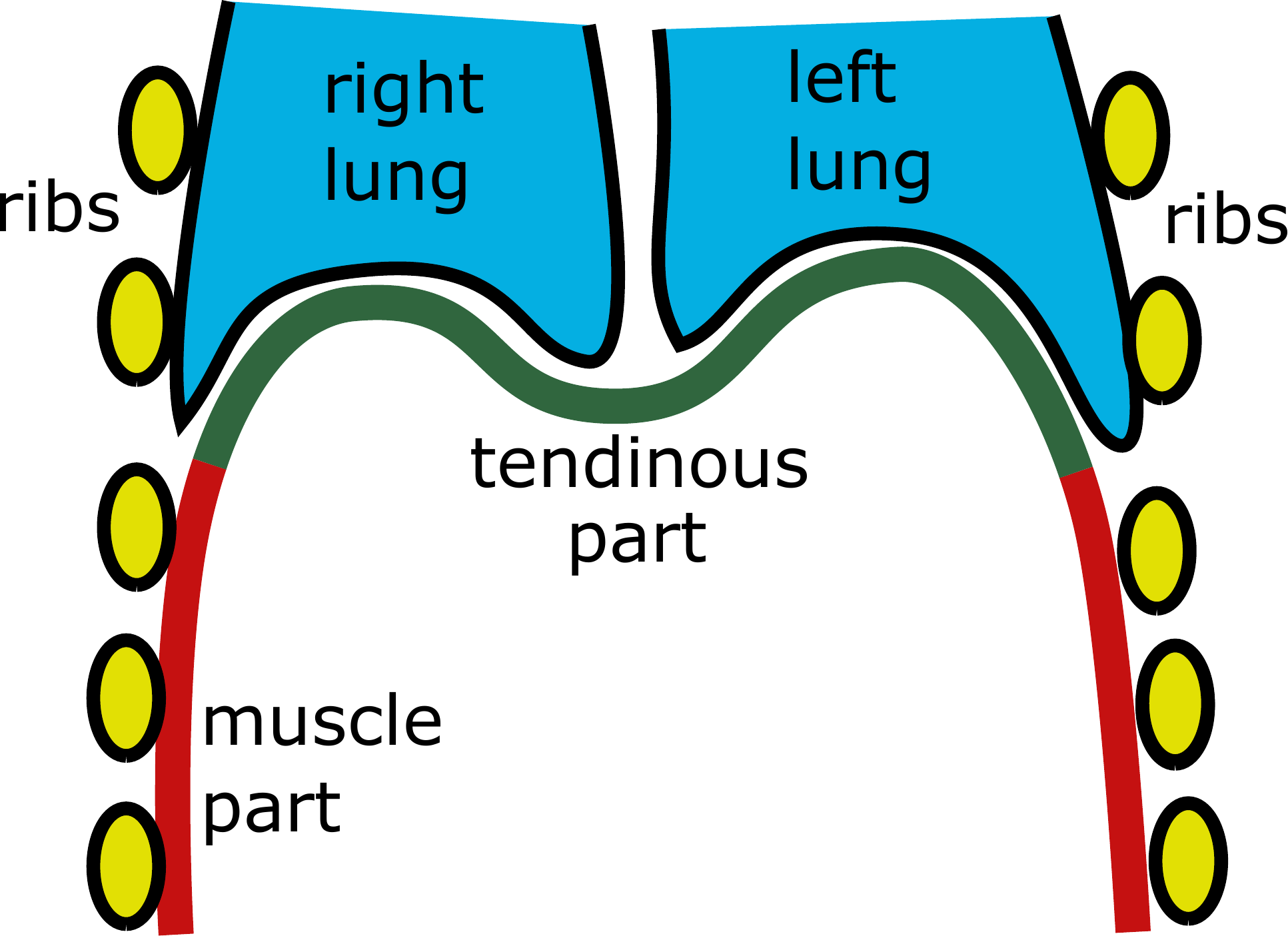}
    \caption{\textbf{Expiration state} with anatomical part}
    \label{fig:1a}
  \end{subfigure}
  \hfill
  \begin{subfigure}[b]{0.48\textwidth}
    \centering
    \includegraphics[height=0.55\textwidth]{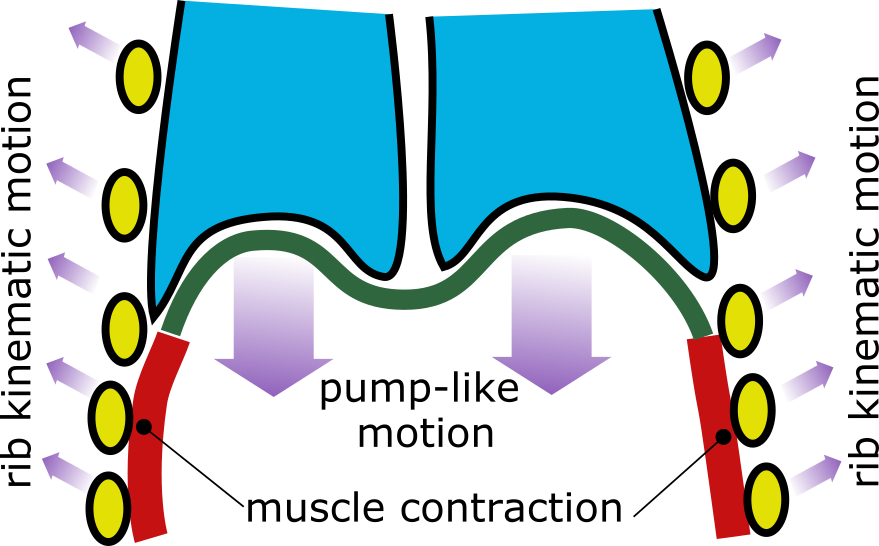}
    \caption{\textbf{Inspiration state} with muscle actions}
    \label{fig:1b}
  \end{subfigure}
  \caption{Anatomy and physiology of the diaphragm (red and green parts) from expiration
to inspiration}
  \label{fig:diaphAnatomy}
\end{figure}

The diaphragm is a musculotendinous structure, approximately double-dome shaped, separating the thoracic and abdominal cavities.  
It is the main muscle of the physiological respiration,  performing 70--80\% of the work of breathing, although it also has non-ventilatory functions, e.g., coughing, hiccups, sneezing, vomiting, and postural functions. It is composed  of three main parts, see Figure~\ref{fig:1a}. 
There are two muscular parts, one median and horizontal in upright position,  separating the thoracic and the abdominal organs, and another lateral muscle 
part, which ends with the costal insertions (attachment to the lower ribs), and there is one central tendon, where the extremities of all the muscle fibers converge. When the diaphragm contracts 
during inspiration (inhalation) with a piston-like motion, the muscle zone thickens (inspiratory thickening), and the domes move caudally (downward) expanding the thorax. Therefore, the air enters the lungs under a pressure gradient, see Figure~\ref{fig:1b}.
%
%

The two-dimensional geometry used for our simulations was extracted from a real patient medical image.
We used a 3D CT scan image (resolution: $0.927 \times 0.927 \times 0.3$ mm$^3$) that was acquired for medical reasons, see our previous work~\cite{invive}. The diaphragm was manually segmented in 3D following visual cues as explained in \cite{villard11}. The labeled voxels were exported to a triangle mesh with the Marching Cube algorithm and the mesh was simplified with a decimation algorithm \cite{levy2010}.

The frontal plane slice we selected is in the middle of the body, corresponding roughly to the anatomic illustration of Figure~\ref{fig:diaphAnatomy}. The raw data consists of a list of 2D vertices where a topology can easily be extracted.
The raw geometry data contains noise from several sources.
There is some CT scan device incertitude (noise, calibration etc) \cite{Zatz77}, the data is the result from a discretized process,  and the labeling data comes from a manual segmentation that is prone to human error. Concerning this last point, the diaphragm is not entirely visible, sometimes part of its voxels also include other organs. The accuracy is linked to the imagination and the anatomical knowledge of the medical expert that performed the segmentation.

\begin{figure}[h!]
\centering
\includegraphics[width=0.85\linewidth]{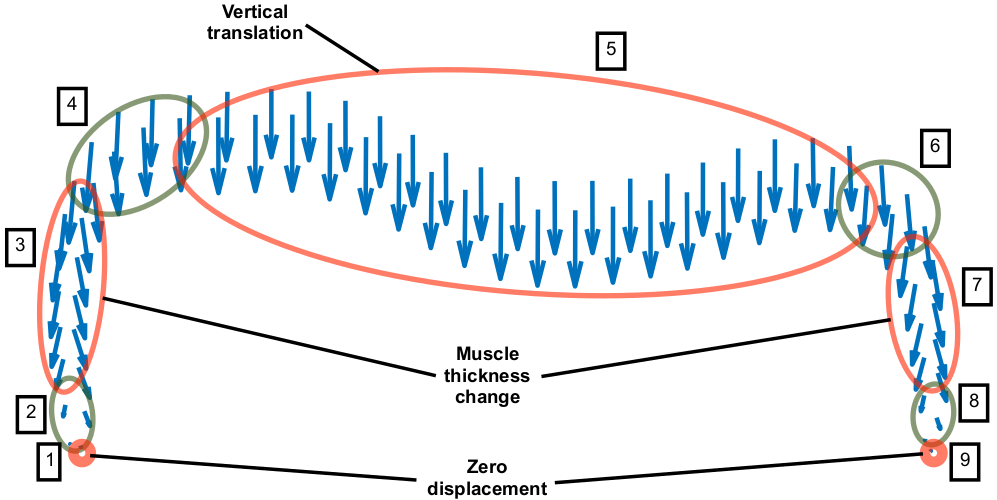}
\caption{The expected displacement in different parts of the slice of the diaphragm. The red areas correspond to the real physiological behavior and the green areas are transition zones. In regions 3 and 7, there is both a thickening of the muscle and a downward motion, while in region 5, the downward motion dominates.}
\label{fig:smoothing:bcPosition}
\end{figure}

%
The expected displacement of the diaphragm between the relaxed and contracted states in the two-dimensional slice is is roughly illustrated in Figure~\ref{fig:smoothing:bcPosition}. For the numerical simulations in Sections~\ref{sec:experiments:dirichlet} and~\ref{sec:experiments:Robin}, we construct boundary conditions to replicate this motion qualitatively.


\section{Equations of linear elasticity}
\label{sec:linearelasticity}
The simplified model that we use is valid for studying 
small deformations of an isotropic and homogenuous diaphragm with a linear elasticity constitutive law. 
The deformation of a diaphragm $\Omega$ is described by applying a displacement field $ u = (u_1, u_2)^T,$ to an object 
location $ y = (y_1, y_2)^T \in \Omega$ such that:
$$ y^* =  y +  u,$$
where $ y^* = (y_1^*, y_2^*)^T \in \Omega^*$ is then a deformed object.
A field derived from the displacements is the stress tensor: 
$$\usigma = \begin{pmatrix}\sigma_{11}& \sigma_{12} \\ \sigma_{21}& \sigma_{22} \end{pmatrix},$$
which measures the internal forces in $\Omega$ as a consequence of deformation. It is related to the displacement field by: 
\begin{equation}
    \label{eq:model:stressAndStrain}
\usigma = \lambda \Tr(\uvarepsilon) + 2\mu\uvarepsilon,\qquad    
    \uvarepsilon = \frac{1}{2}\left( \left(\nabla u\right)^T + \nabla  u \right),
\end{equation}
where $\uvarepsilon\in\mathbb{R}^{2\times 2}$ is the strain tensor and $\Tr(\uvarepsilon) = \left(\varepsilon_{11} + \varepsilon_{22}\right)I$ is its trace. 
The scalars $\lambda$ and $\mu$ are the Lam\'e parameters, which are related to the Young modulus $E$ and the Poisson ratio $\nu$ of the material through:
\begin{equation}
    \label{eq:model:lame}
    \lambda = \frac{E\nu}{\left(1+\nu\right)\left(1-2\nu\right)},\qquad 
    \mu = \frac{E}{2\left(1+\nu\right)}.
\end{equation}
For our computations 
we use $E=10^5\, \text{Pa}$ and $\nu = 0.3$.
A special stress measure that is also of interest to us is the Von Mises stress:
\begin{equation}
    \label{eq:model:vonMises}
\sqrt{\sigma_{11}^2 - \sigma_{11}\sigma_{22} + \sigma_{22}^2 +  3\sigma_{12}^2},
\end{equation}
which provides a scalar measure of the total stress.

The equations of linear elasticity are derived from the force equilibrium (Newton's second law) imposed on $\Omega$. 
We have:
\begin{equation}
\begin{aligned}
\label{eq:model:secondNewton}
- \nabla \cdot \usigma =  f \text{ on } \Omega,
\end{aligned}
\end{equation}
where $ f = \left(f_1( y), f_2( y)\right)^T$ is a field of 
internal forces in the horizontal and the vertical direction.
The Dirichlet and traction boundary conditions are prescribed on two disjoint parts of the 
boundary, which together form the whole boundary $\partial\Omega = \partial\Omega_0 \cup \partial\Omega_1$:
\begin{equation}
    \begin{aligned}
\label{eq:model:secondNewton_bcs}
 u &=  g \text{ on } \partial\Omega_0, \\
\usigma \cdot  n &=  h \text{ on } \partial\Omega_1.
    \end{aligned}
\end{equation}
The first boundary condition with the right hand side $ g = \left(g_1( y), g_2( y)\right)^T$ is the displacement 
condition. 
The second boundary condition with the right hand side $ h = \left(h_1( y), h_2( y)\right)^T$ is the traction condition. 
An equivalent form of \victor{\eqref{eq:model:secondNewton_bcs}} is the Robin boundary condition:
\begin{equation}
\begin{aligned}
    \label{eq:model:secondNewton_bcs_Robin}
     u\, \kappa_0(y)  +  \left(\usigma \cdot  n\right)\, \kappa_1(y) =   g\, \kappa_0(y) +   h\, \kappa_1(y)  \text{ on } \partial\Omega = \partial\Omega_0 \cup \partial\Omega_1,
\end{aligned}
\end{equation}
where $\kappa_0$ and $\kappa_1$ correspond to two spatially dependent coefficients, in this case discontinuous: 
$\kappa_0(y) = 1$ when $y\in \victor{\partial} \Omega_0$ and zero otherwise, and $\kappa_1(y) = 1$ when $y\in \victor{\partial}\Omega_1$ and zero otherwise.

\victor{In Section \ref{sec:experiments:dirichlet}} we first compute solutions using a Dirichlet condition for the whole boundary. 
($\kappa_0(y) = 1$, $\kappa_1(y)=0$, $y \in \partial \Omega$). 
Then, \victor{in Section \ref{sec:experiments:Robin}}, we compute 
solutions using smoothed Dirichlet and traction conditions (smooth Robin coefficients), which approximates the 
imposition of these two conditions on two disjoint parts of the diaphragm. However, since physiology suggests a smooth solution, 
we see this as modeling rather than as an error.

For all of the computations in this paper we use the displacement formulation of \eqref{eq:model:secondNewton}, obtained by 
using the relation between stress and displacement  given in \eqref{eq:model:stressAndStrain}. The force equilibrium 
\eqref{eq:model:secondNewton} then expands to:
$$- \nabla \cdot \usigma = - \mu \nabla^2 {u} - (\lambda + \mu)\nabla (\nabla \cdot {u}) =  f,$$
and the traction boundary condition from \eqref{eq:model:secondNewton_bcs} to:
$$\usigma \cdot  n = \left[\lambda \left(\nabla \cdot  u \right)I + \mu \left(\left(\nabla  u\right)^T + \nabla  u\right)\right] \cdot  n =  h.$$
The displacement formulation of the linear elasticity equations with the boundary condition from \eqref{eq:model:secondNewton_bcs_Robin} 
is then written as:
\begin{equation}
\begin{aligned}
    \label{eq:model:naviercauchy}
- \mu \nabla^2 u - (\lambda + \mu)\nabla (\nabla \cdot u) = &\, f  &\text{ on } \Omega,& \\
u\, \kappa_0 + \left(\left[\lambda \left(\nabla \cdot  u \right)I + \mu \left(\left(\nabla  u\right)^T + \nabla  u\right)\right] \cdot  n \right)\kappa_1 = &\, g\, \kappa_0 + h\, \kappa_1 & \text{ on } \partial\Omega&.
\end{aligned}
\end{equation}
For simplicity we rewrite the system above as:
$$Du(y) = F(y),$$ 
where:
\begin{equation}
\begin{aligned}
    \label{eq:model:PDE}
    D u( y) &=
    \left\{
    \begin{array}{ll}
        D_2  u( y), &  y \in \Omega \\
        \kappa_0( y) D_0  u( y) + \kappa_1( y) D_1  u( y) \cdot  n( y), &  y \in \partial\Omega
    \end{array}
    \right. & 
    \\
    F(y) &=
    \left\{
    \begin{array}{ll}
         f( y), &  y \in \Omega  \\
         \kappa_0( y) g( y) + \kappa_1( y) h( y), &  y \in \partial\Omega
    \end{array}
    \right. &
\end{aligned}
\end{equation}
Here $D_2$, $D_1$, $D_0$ are the expanded operators that in \eqref{eq:model:PDE} correspond to 
$\Omega$, $\partial\Omega_1$ and $\partial\Omega_0$ respectively. If we let $\nabla_{ij}=\frac{\partial^2}{\partial y_i\partial y_j}$ and  $\nabla_i=\frac{\partial}{\partial y_i}$, then the operator $D_2$ is:
\begin{equation}
\label{eq:model:PDE_D_2}
D_2 
= 
-
\begin{pmatrix}
(\lambda+2\mu)\nabla_{11} + \mu \nabla_{22} & (\lambda + \mu) \nabla_{12} \\
(\lambda + \mu) \nabla_{12} & \mu \nabla_{11} + (\lambda+2\mu) \nabla_{22}
\end{pmatrix},
\end{equation}
and the boundary operators $D_1$ and $D_0$ are:
\begin{equation}
\label{eq:model:PDE_D_1_0}
D_1 
= 
\begin{pmatrix}
\mu n_2 \nabla_2 + \left(\lambda + 2\mu\right) n_1 \nabla_1 & \lambda n_1 \nabla_2 + \mu n_2 \nabla_1 \\
\mu n_1 \nabla_2 + \lambda n_2 \nabla_1 & (\lambda+2\mu)n_2\nabla_2 + \mu n_1 \nabla_1
\end{pmatrix}, \quad 
D_0 
=
\begin{pmatrix}
I & 0 \\ 0 & I
\end{pmatrix}.
\end{equation}

\section{Approximation of functions using the RBF-FD method}
\label{sec:rbffd}
Here we describe the approximation framework which is used when representing the unknown displacement and the 
stress (which are the outputs of the simulation), and the known data vectors (the input data to the simulation) as continuous functions. 
The following discussion is general in the sense that 
the data vector can be known or unknown. Concrete examples of using the framework developed in this section are given 
in later sections.

Given a data vector $\tilde u(X) = [\tilde u(x_1), \tilde u(x_2), ..., \tilde u(x_N)]^T$ where $x_i \in \mathbb{R}^d$ 
we construct a semi-discrete evaluation operator $E=E(y,X)$, which, for 
any point $y \in \Omega \subset \mathbb{R}^d$, returns the value of the interpolant of the data $\tilde u(X)$ at that point, 
such that:
\begin{equation}
    \label{eq:rbffd:evaluation}
\tilde u(Y) = E(y,X) \tilde u(X).
\end{equation}
Now we use \eqref{eq:rbffd:evaluation} for each $y \in Y$ to obtain a system of equations and form the matrix 
$E(Y,X)$ of size $M \times N$, $M \geq N$. Note that in our notation this is equivalent to setting $y=Y$.
The RBF-FD procedure constructs the matrix $E(Y,X)$ using a sequence of local, stencil-based interpolation problems 
which are exact for a cubic or quintic polyharmonic spline basis (PHS) and a (multivariate) monomial basis of degree $p$.
The matrix $E(Y,X)$ is sparse, that is, 
it contains $n \ll N$ non-zero elements per row, where $n$ is the stencil size defined by:
\begin{equation}
\label{eq:rbffd:stencilsize}
n = 2\,\binom{p+d}{d}.
\end{equation}
Details about using PHS plus the monomial basis in RBF-FD approximations can be found in \cite{BFFB17, BFF19, Bayona19, Barnett15}.
\victor{In the same way,} we construct a semi-discrete oprator for differentiation $D^{\mathcal{L}} = D^{\mathcal{L}}(y,X)$ which locally evaluates any 
derivative $\mathcal{L}$ of $\tilde u(Y)$:
\begin{equation}
    \label{eq:rbffd:differentiation}
    \mathcal{L} \tilde u(Y) = D^{\mathcal{L}}(y,X) \tilde u(X).
\end{equation}

\victor{Both matrices, $E(Y,X)$ and $D^{\mathcal{L}}(Y,X)$, can be formed using the MATLAB code available in \cite{tominec_rbffdcode}.} 
\victor{For completeness, we below provide the steps to compute the matrix elements (weights) and to assemble 
$D^{\mathcal{L}}(Y,X)$ and $E(Y,X)$ 
for the point sets $Y$ and $X$ in any dimension $d$.}
\begin{enumerate}
    \item Let $x_1^{(k)} = x_i \in X$ (one stencil center) and find 
    $n$ closest neighbors $\left\{x^{(k)}_j\right\}_{j=1}^n$ (stencil points) around it using the Euclidean distance.
    \item Scale and shift the stencil points to a unit domain $[-1,1]^d$. \victor{Save the scaling as $s^{(k)}$.}
    \item Form a square interpolation matrix $A^{(k)}$, where $A^{(k)}_{ij}=r^3=||x_i^{(k)}-x_j^{(k)}||_2^3$ 
    \victor{and $x_i^{(k)},x_j^{(k)},\, i,j = 1,..,n$ belong to the stencil}.
    \item Form a rectangular polynomial matrix $P^{(k)}$, where $P^{(k)}_{il} = p_l(x_i^{(k)})$, $i=1,..,n$, $l=1,..,m$ is a sampled $d$-dimensional monomial basis 
    with $m$ basis functions.
    \victor{When $n$ is chosen as in \eqref{eq:rbffd:stencilsize} 
    then $m = \frac{n}{2}$ and the size of the matrix $P^{(k)}$ is $n \times \frac{n}{2}$.}
    \item Using $A^{(k)}$ and $P^{(k)}$, form the augmented local interpolation matrix:
    \begin{equation}
	    \label{eq:M}
	    \tilde{A}^{(k)} = \begin{pmatrix}
		    A^{(k)} & P^{(k)} \\
	    	(P^{(k)})^{T} & 0
	    \end{pmatrix}
    \end{equation}
    
    \item Repeat steps 1--5 for every $x \in X$ in order to form all $\tilde{A}^{(k)}$.
\end{enumerate}
The local evaluation and differentiation weights can then be computed in the following way:
\begin{enumerate}    
    \item Take one evaluation point $y_l \in Y$ and find the first closest point from the $X$ point set. We denote it by $x_1^{(k)}$.
    \item Scale $y_l$ to a unit domain $[-1,1]^d$ \victor{using the previously computed scaling $s^{(k)}$}.
    \item Form a vector $b_1 = \mathcal{L}||y_l-x^{(k)}_j||_2^3$, where $\{x^{(k)}_j\}_{j=1}^{n}$ is the local neighborhood of the 
    center point $x_1^{(k)}$, \victor{where $\mathcal{L} = 1$ for constructing evaluation weights, or a derivative operator for constructing differentiation weights.}
    \item Form a vector $b_2 = \mathcal{L} p_j(y_l),\, j=1,..,\frac{n}{2}$.
    \item Concatenate the two vectors into $b(y_l) = [b_1,\, b_2]$.
    \item Use the augmented interpolation matrix that belongs to $x_1^{(k)}$ and compute the local weights by using the relation 
    $w_{\mathcal{L}}(y_l) = (\tilde{A}^{(k)})^{-1} b_{\mathcal{L}}(y_l)$.
    \item Store $w_{\mathcal{L}}(y_l)$ in the $l$-th row of the matrix $W_{\mathcal{L}}(Y)$.
    \item Repeat steps 1--6 for every $y_l \in Y$ in order to form all local weights.
\end{enumerate}
\victor{When using these steps to compute the weights, it is possible to avoid storing the matrix 
$\tilde{A}^{(k)}$ by combining the two parts such that the evaluation/differentiation weights are computed for all $y_l$ 
which select $x_1^{(k)}$ as its closest stencil center, immediately after step 5 of the first list.}
Once the matrix of local weights $W_{\mathcal{L}}(Y)$ of size $M \times n$ is computed, 
these weights are assembled into a global rectangular matrix, $E(Y,X)$ for evaluation and 
$D_{\mathcal{L}}(Y,X)$ for differentiation, both of size $M \times N$. This can be done by using a matrix $\Gamma \in \mathbb{Z}^{N\times n}$, 
which is a list of indices of the local neighborhoods of points around every stencil point $x_1^{(k)} \in X$. 
Additionally, a matrix $\kappa \in \mathbb{Z}^{M\times 1}$ is needed, which is a list of indices of the 
stencil centers that are closest to each evaluation point $y_l$. $\Gamma$ and 
$\kappa$ are found using the k-nearest neighbor method with $k$ equal to $n$ and $1$, respectively.
Using the MATLAB programming language, the two lists and the sparse matrix can be efficiently computed by 
invoking the following three commands:
\begin{verbatim}
    Gamma = knnsearch(X,X,'k',n)
    kappa = knnsearch(X,Y,'k',1)
    D_L = sparse(repmat(1:M, 1, n), Gamma(kappa,:), W_L, M, N, N*n)
\end{verbatim}
where the first two arguments to the function \texttt{sparse()} have the same shape as $W_{\mathcal{L}}(Y)$ 
and contain the row and column indices for inserting the locally computed weights into the global matrix.

\section{Discretization of the linear elastic model using the uniftted RBF-FD method}
\label{sec:unfitted_discretization}
In this section we employ the unfitted RBF-FD method \cite{tominec2020unfitted} to discretize the linear elasticity equations 
\eqref{eq:model:PDE}, \eqref{eq:model:PDE_D_2} and \eqref{eq:model:PDE_D_1_0} over the diaphragm. 
The method relies on constructing rectangular 
differentiation matrices $D^\mathcal{L}(Y,X)$ as described in 
Section \ref{sec:rbffd}. The matrices then replace the differential operators in the PDE, together with the corresponding 
discrete right hand sides $F(Y)$.

\begin{figure}[h!]
    \centering
\includegraphics[width=0.4\linewidth]{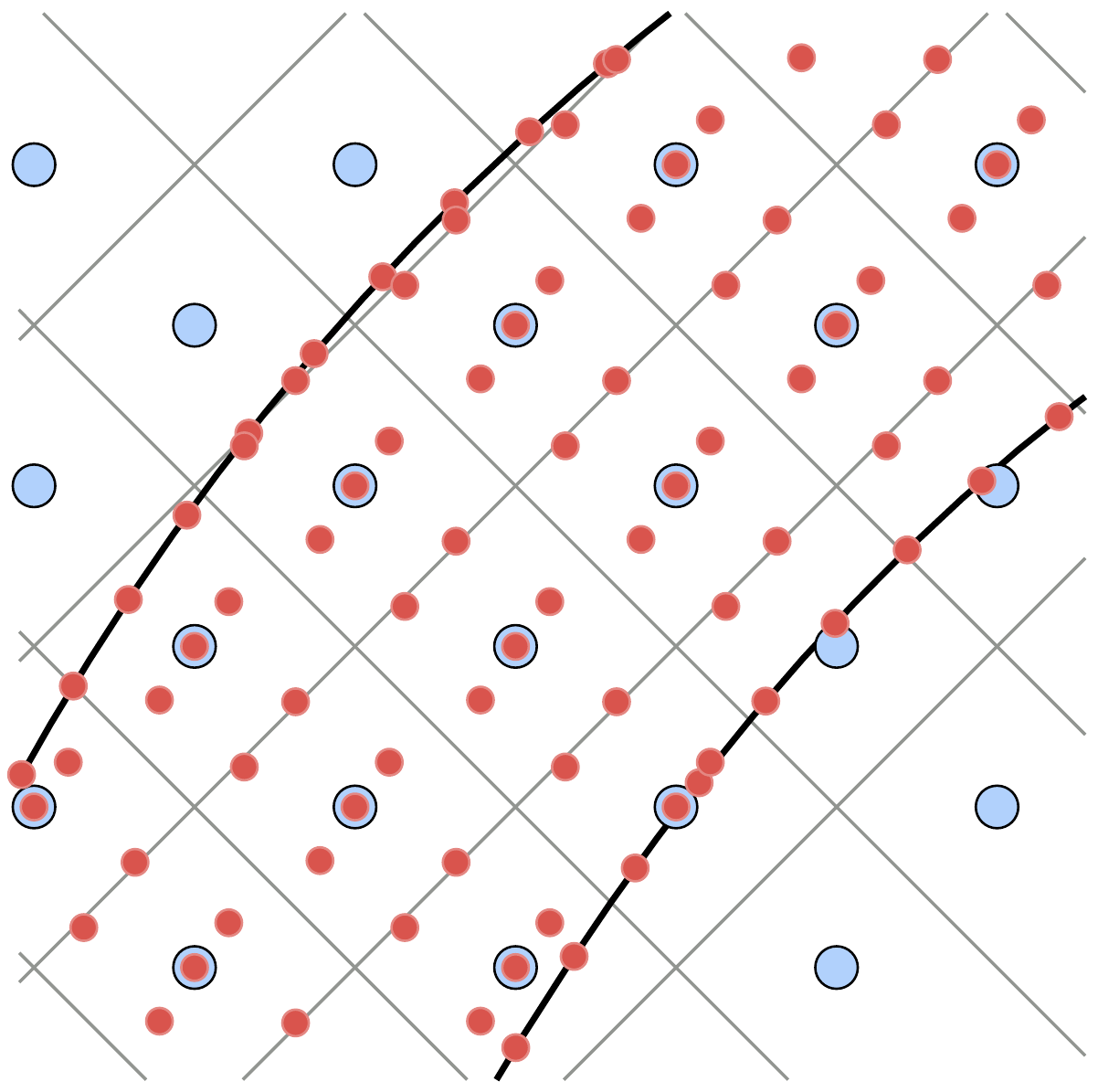}
    \caption{The black curve represents a part of the boundary of the diaphragm. 
    Node points in $X$ (blue markers) and the corresponding Voronoi regions (grey lines) are shown 
    together with interior evaluation points, $Y_1$, and boundary evaluation points, $Y_b$, (red markers). 
    The same template of interior evaluation points is used in each Voronoi region inside the diaphragm geometry.}
    \label{fig:method:diaphragm_evalNodes}
\end{figure}

To form the differentiation matrices, we first decide the degree $p$ of the monomial basis that we are appending to the 
PHS approximation, and compute the stencil size $n$ using \eqref{eq:rbffd:stencilsize}.
\begin{figure}[h!]
    \centering
    \includegraphics[width=0.8\linewidth]{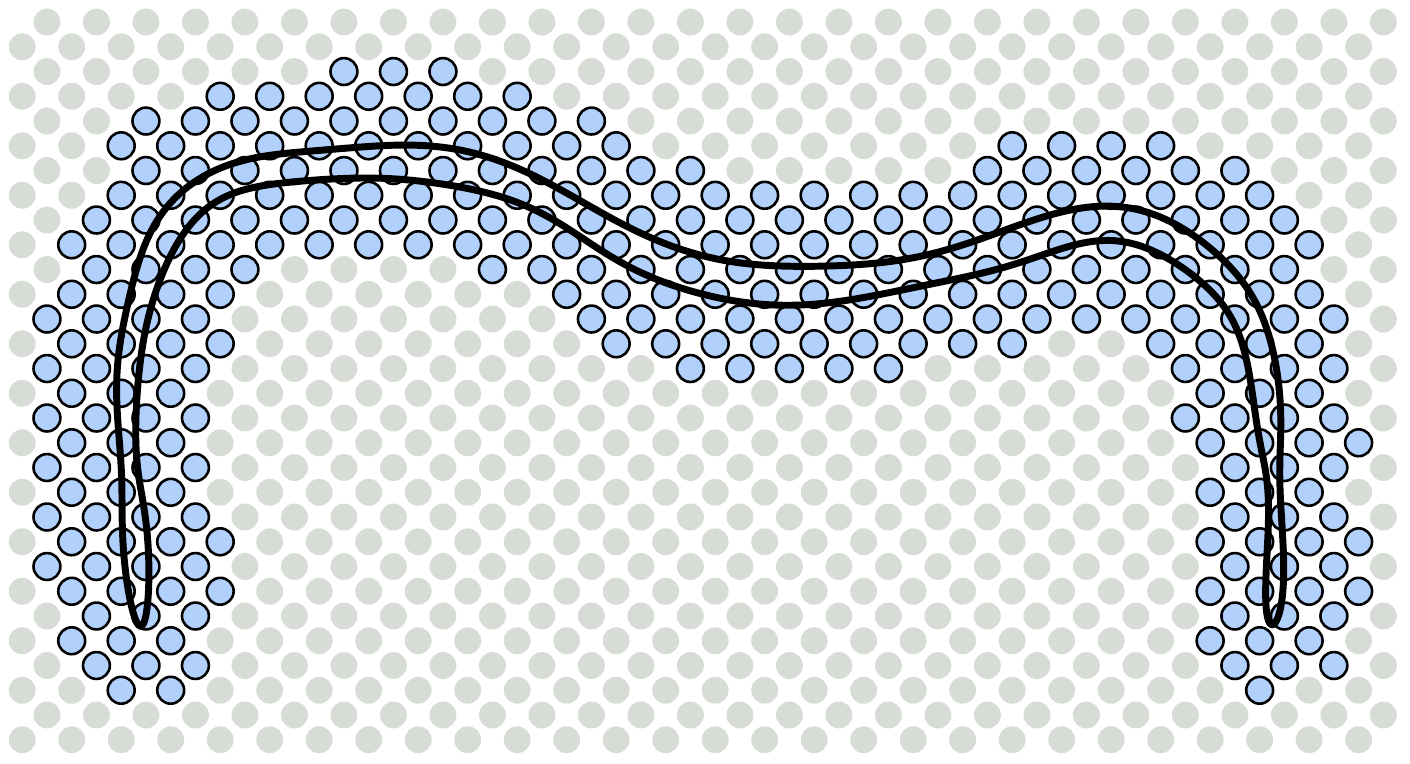} 
    \caption{The initial point set $X_1$ (tilted Cartesian points) is distributed over a box that encloses the boundary 
    of the diaphragm (black curve). Points that are more than half a stencil size away from the geometry (grey) are removed. 
    The remaining points (blue) form the point set $X$.}
    \label{fig:method:diaphragm_interpNodes}
\end{figure}
Then we construct an interpolation point set $X$ 
that extends over the diaphragm (see Figure \ref{fig:method:diaphragm_interpNodes}), and an 
evaluation point set $Y$ (see Figure \ref{fig:method:diaphragm_evalNodes}) that conforms to the geometry of the diaphragm. 
Those two point sets are obtained in four simple steps:
\begin{itemize}
\item The initial point set $X_1$ is a tilted Cartesian node layout with spacing $h$ in a box that 
encloses the diaphragm. \victor{(the grey and blue points in Figure \ref{fig:method:diaphragm_interpNodes})}.
\item Then the point set $Y_1$ is generated by placing $q$ points with average spacing $h_y$ in each Voronoi region inside the diaphragm, 
defined by the points in $X_1$. \victor{(the red points in Figure \ref{fig:method:diaphragm_evalNodes})}.
\item In addition, the point set $Y_b$ with the same average spacing $h_y$ is generated by placing points 
along the boundary of the diaphragm, see Figure \ref{fig:method:diaphragm_evalNodes}.
\item The final evaluation point set is given by $Y = Y_1 \cup Y_b$, see Figure \ref{fig:method:diaphragm_evalNodes}.
\item The final node set $X$ (the blue points in Figure \ref{fig:method:diaphragm_interpNodes} and Figure \ref{fig:method:diaphragm_evalNodes}) 
is formed by reducing $X_1$ by removing points that fall more than half stencil size outside the diaphragm. In Matlab this can be done by:
\begin{verbatim}
    X = X_1(unique(knnsearch(X_1,Y,’k’,ceil(0.5*n)), :);
\end{verbatim}
\end{itemize} 
The last step ensures that the columns of $E(Y,X)$ and $D^{\mathcal{L}(Y,X)}$ are non-zero \cite{tominec2020unfitted}. 
Note that the evaluation point set does not need to be constructed by placing precisely $q$ points in every Voronoi region. 
It is possible to use any 
global point set, which is quasi-uniform by nature (e.g. Halton points), 
and as such on average samples every Voronoi region with approximately $q$ points.


\elisabeth{%
Next, we use the method described in Section~\ref{sec:rbffd} to discretize the continuous operator $D_2$ in~\eqref{eq:model:PDE_D_2} using the interior evaluation points $Y_1$, and to discretize the continuous operators $D_1$ and $D_0$ in~\eqref{eq:model:PDE_D_1_0} using the boundary evaluation points $Y_b$. We let $D_k^{ij}(\,\cdot\,,X)$ denote the differentiation matrix that approximates element $i,j$, of the operator $D_k$, $k=0,\,1,\,2$. If we let $i,j=1,2$, while $i\neq j$, then we can express the differentiation matrices as:
  \begin{eqnarray*}
    D_2^{ii}(Y_1,X) &=& (\lambda+2\mu)D^{\nabla_{ii}}(Y_1,X) + \mu D^{\nabla_{jj}}(Y_1,X), \\
    D_2^{ij}(Y_1,X) &=& (\lambda + \mu) D^{\nabla_{ij}}(Y_1,X), \\
    D_0^{ii}(Y_b,X) &=& E(Y_b,X), \\
    D_0^{ij}(Y_b,X) &=& 0, \\
    D_1^{ii}(Y_b,X) &=& (\lambda+2\mu)n_i D^{\nabla_i}(Y_b,X) + \mu n_j D^{\nabla_j}(Y_b,X),\\
    D_1^{ij}(Y_b,X) &=& \lambda n_i D^{\nabla_j}(Y_b,X) + \mu n_j D^{\nabla_i}(Y_b,X). \\
\end{eqnarray*}
}%
\elisabeth{We express the discrete Robin coefficients as $K_i(Y_b) = \diag(\kappa_i(Y_b))$, $i=0,\,1$, and introduce a scaling $\beta_i$ for equations connected with the operator $D_i$.}
\elisabeth{Finally, we form the rectangular system of size $2M \times 2N$ that discretizes  the Navier-Cauchy system with 
 Robin boundary conditions~\eqref{eq:model:PDE}:}
 \begin{flalign}
   \small
   \label{eq:method:DandF}
    \begin{pmatrix}
        \beta_2\,D_2^{11} & \beta_2\, D_2^{12} \\
        \beta_2\, D_2^{21} & \beta_2\, D_2^{22} \\
        \beta_0\,K_0\, D_0^{11}+  \beta_1\,K_1\, D_1^{11}\hspace*{-2mm}& \beta_1\,K_1\, D_1^{12}\\    
        \beta_1\,K_1\, D_1^{21} & \hspace*{-2mm}\beta_0\,K_0\, D_0^{22} + \beta_1\,K_1\, D_1^{12} 
    \end{pmatrix}
    \begin{pmatrix}
    \tilde u_1(X)\\
    \tilde u_2(X)
    \end{pmatrix}
    =
    \begin{pmatrix}
        \beta_2\, f_1(Y_1)\\
        \beta_2\, f_2(Y_1)\\
        \beta_0\, K_0\, g_1(Y_b) + \beta_1\, K_1\, h_1(Y_b)\\
        \beta_0\, K_0\, g_2(Y_b) + \beta_1\, K_1\, h_2(Y_b)\\
    \end{pmatrix}.
\end{flalign}
 In the numerical experiments, we use the following scale factors:
\begin{equation}
    \beta_2 = \frac{1}{\mu}\, h_y,\quad \beta_1 = \frac{10}{\mu}\,\frac{1}{h}\, h_y^{\frac{1}{2}}, \quad \beta_0 = \frac{1}{h}\, h_y^{\frac{1}{2}},
\end{equation}
\victor{where $h$ and $h_y$ are the average internodal distances in the $X$ and $Y$ point sets, respectively. These are computed as:
\begin{equation}
    \label{eq:method:h_and_hy}
h = \frac{1}{N} \sum_{j=1}^N \min_{i \neq j} \|x_i - x_j\|_2,\quad h_y = \frac{1}{M} \sum_{j=1}^M \min_{i \neq j} \|y_i - y_j\|_2.
\end{equation}
}%
\new{The choice of scale factors is based on the papers \cite{tominec2020unfitted, ToLaHe20}, where the scale factors $h_y$ for the interior and $h_y^{1/2}$ for the boundary are used, such that the norms of the discrete least squares problem approximate the continuous $L_2$-norm. The additional $1/h$ scaling increases the weight of the boundary conditions and improves convergence. The Lam\'e parameter $\mu$ is large and affects the scaling between different equations. Therefore, we also include the factor $10/\mu$ in $\beta_1$.
}

We solve the least-squares problem \eqref{eq:method:DandF}, for the nodal values 
$\tilde u_1(X)$ and $\tilde u_2(X)$, using backslash in MATLAB. 
After that the solution is evaluated at the $Y$ point set, using the evaluation matrix $E$:
\begin{equation}
    \tilde u_1(Y) = E(Y,X)\, \tilde u_1(X),\, \qquad \tilde u_2(Y) = E(Y,X)\, \tilde u_2(X).
\end{equation}
The strains and the stresses \eqref{eq:model:stressAndStrain}, \eqref{eq:model:vonMises} are computed by applying appropriate 
differentiation matrices to the solution coefficients $\tilde u_1(X)$ and $\tilde u_2(X)$.

\section{Smoothing of geometry and boundary data}
\label{sec:smooth}
If we view the two-dimensional diaphragm from the continuous perspective, the geometry can be described as a closed curve. We choose to parametrize this curve by $t\in[0,\,2\pi]$, where the starting point $t=0$ is the same as the final point $t=2\pi$. To benefit from the potential high-order convergence of the unfitted RBF-FD method, we need to approximate the boundary curve and the boundary data as functions of $t$ with enough smoothness that the convergence of the PDE problem is not adversely affected.
%
To avoid reduced accuracy due to boundary errors near the artificial end points of the interval, we extend the domain periodically to $t\in[-2\pi,\,4\pi]$ for the approximation.
We discretize the extended domain using the uniformly spaced node points $T=\{t_j\}_{j=1}^{N_g}$.
Given $\tilde{M}_g$ data points, we replicate these periodically to get the extended data set $(T^d,G)=\{(t^d_i,g_i)\}_{i=1}^{M_g}$, where $M_g=3\tilde{M}_g>N_g$. We use a one-dimensional RBF-FD approximation, based on a quintic PHS basis augmented with polynomials of degree $p_g=6$,
to form an overdetermined linear system for the nodal values $g(T)$
\begin{equation}
  E_g(T^d,T)g(T)=G.
\end{equation}
%
%
To enforce continuity at $t=0$, we add the following equality constraints
\begin{equation}
  \frac{d^sg(0)}{dt^s} - \frac{d^sg(2\pi)}{dt^s} = 0,\quad s=0,\ldots,p_g-1,
  \label{eq:constr}
\end{equation}  
and solve the constrained least squares problem
\begin{equation}
    \label{eq:boundary:system}
\begin{pmatrix}
2 E_g^T E_g & B_g^T \\
B_g & 0
\end{pmatrix}
\begin{pmatrix}
    g(T) \\
\lambda
\end{pmatrix}  
=  
\begin{pmatrix}
    2 E_g^T G \\
    0
\end{pmatrix},
\end{equation}
where $B_g$ contains the $p_g$ constraints~\eqref{eq:constr} and $\lambda$ contains the corresponding Lagrange multipliers. 

\begin{figure}[h!]
    \centering
    \hspace*{-1mm}\includegraphics[width=0.51\linewidth]{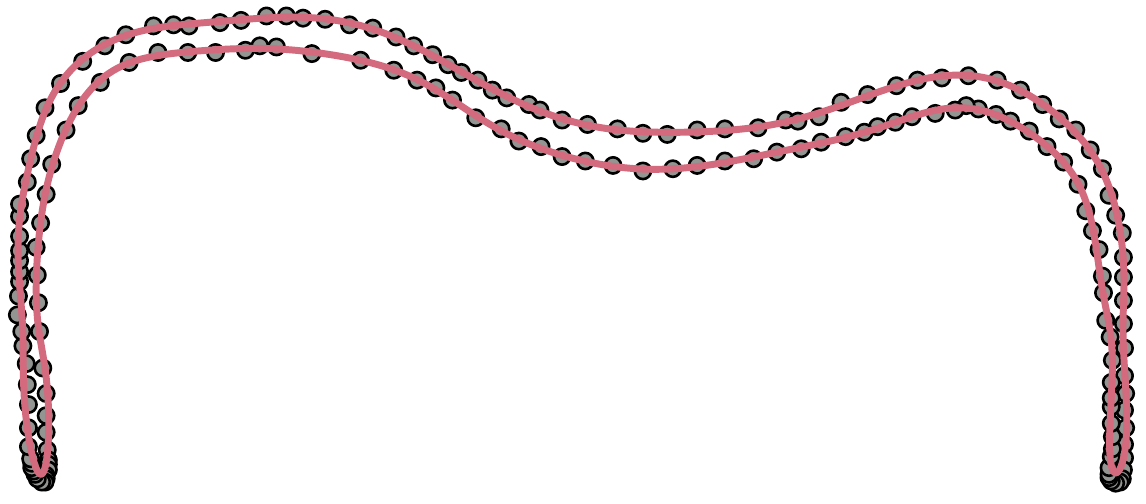} \\
    \includegraphics[width=0.65\linewidth]{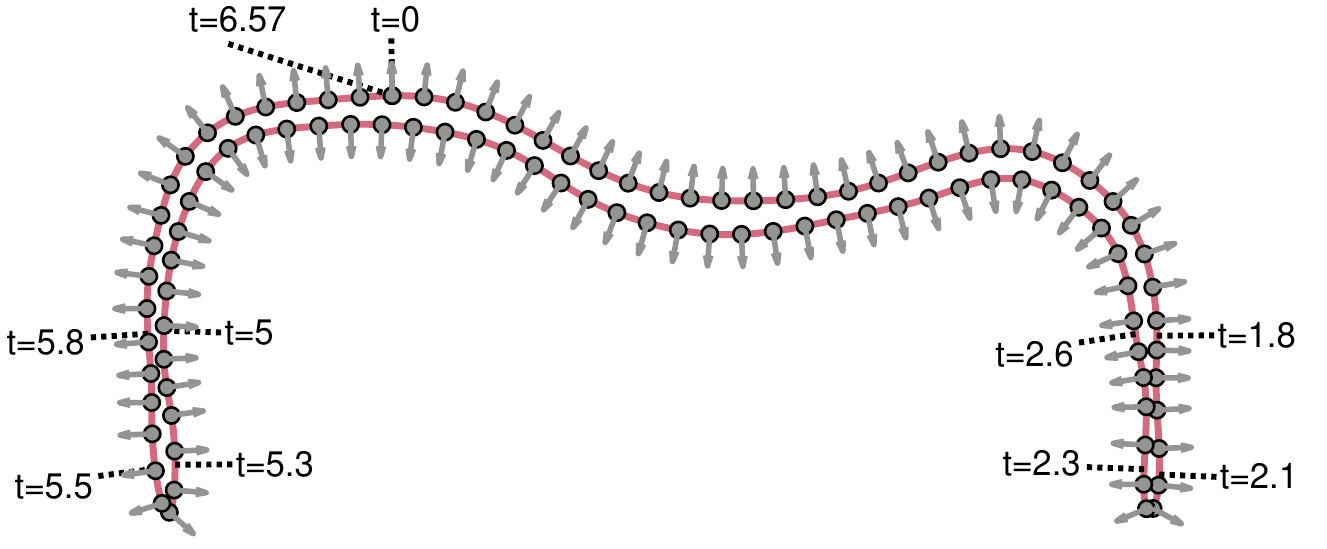}
    \caption{The smoothed boundary geometry curve (solid line) is shown in both subfigures. 
    The curve was computed using $N_g=133$ node points, $M_g=177\cdot 3=531$ data points, and stencil size $n=28$. 
    The markers show the $\tilde{M}_g=177$ vertex data points (top) and uniform evaluation points (bottom). 
    The normals computed from the approximation, as well as the values of the parameter $t$ along the curve, are also shown in the bottom subfigure.
    }
    \label{fig:boundary}
\end{figure}

\begin{figure}[h!]
    \centering
    \includegraphics[width=0.49\linewidth]{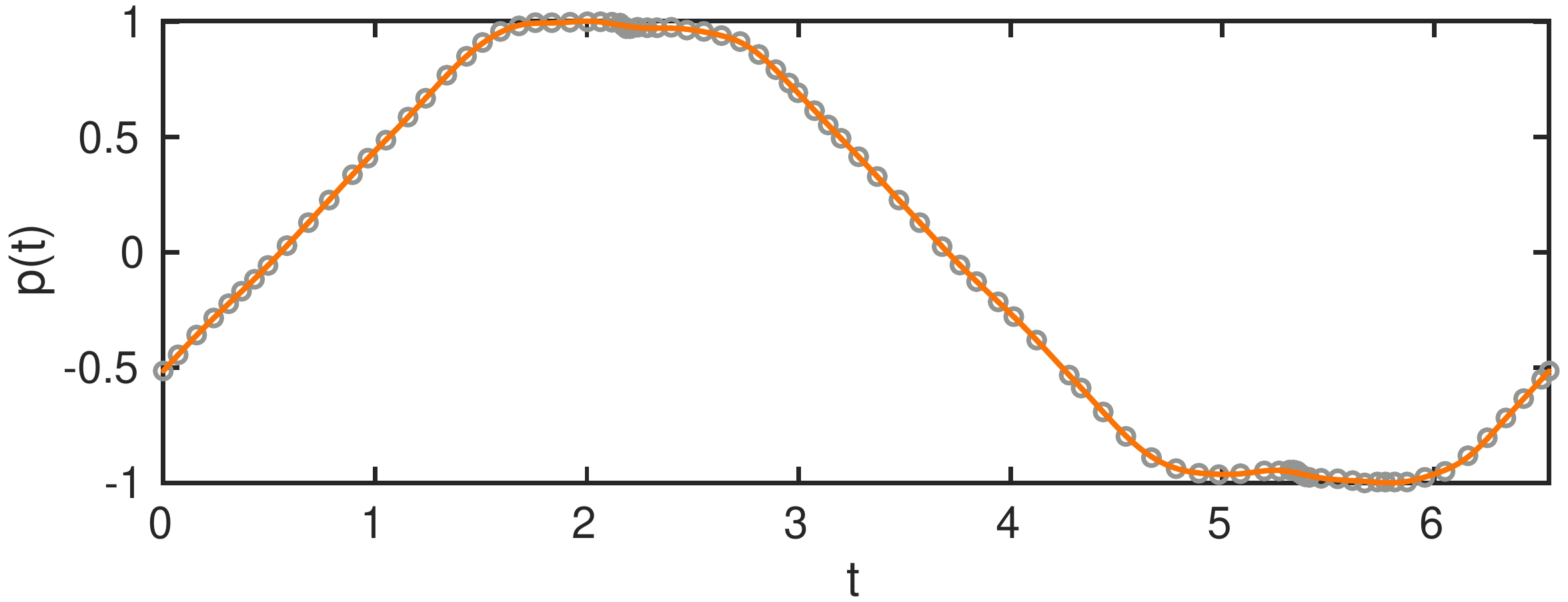}
    \includegraphics[width=0.49\linewidth]{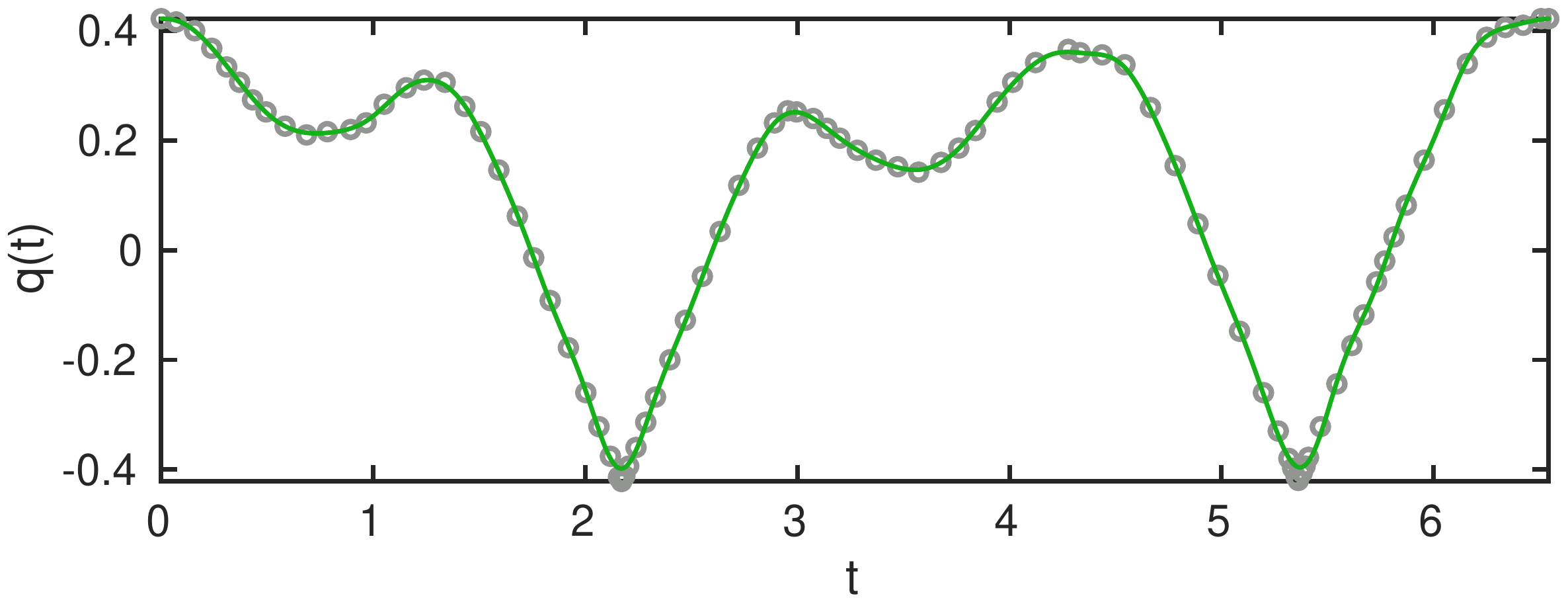}
    \caption{The two smoothed coordinate functions approximating the geometry. The function $p(t)$ (left) corresponds to the horizontal coordinate and $q(t)$ (right) corresponds to the vertical coordinate. The markers show the initial data locations.} 
    \label{fig:p_and_q}
\end{figure}

To find the smooth boundary curve from the initial vertex data $\tilde{x}_i^d$, $i=1,\ldots,\tilde{M}_g$, we first scale the data such that $x_i^d=(p_i,q_i)=s_\Omega\tilde{x}_i^d$, $i=1,\ldots,\tilde{M}_g$, where $s_\Omega=156.92^{-1}$\text{mm}$^{-1}$. 
The scaling was chosen such that all data points fall within $[-1,\,1]^2$. Then we compute an approximate arclength parametrization using the Euclidean distance between the scaled vertices, such that $t_i^d=2\pi\sum_{j=1}^{i-1}\|x_{j+1}^d-x_{j}^d\|/\sum_{j=1}^{\tilde{M}_g}\|x_{j+1}^d-x_j^d\|$, where $x_{\tilde{M}_g+1}^d=x_1^d$. Then we replicate the data over the extended domain.
Finally, system~\eqref{eq:boundary:system} is solved for each coordinate function $p(t)$ and $q(t)$. The resulting boundary curve is shown in Figure~\ref{fig:boundary} and the individual coordinate functions are shown in Figure~\ref{fig:p_and_q}. Since the curve parametrization here is in the clockwise direction, 
the outward normals are computed as $n(t) = (-q'(t),p'(t))/\|(-q'(t),p'(t))\|$. 
%
%

For the boundary data functions $g$ and $h$ in~\eqref{eq:model:naviercauchy}, we manufacture data to mimic 
the expected physiological behaviour shown in Figure~\ref{fig:smoothing:bcPosition}. A few data points are placed in the regions where we have some 
information (regions 1, 3, 5, 7, and 9 in Figure~\ref{fig:smoothing:bcPosition}), and then the rest of the data is generated through linear interpolation. 
This results in gradual transitions in the regions where we lack information. The data points and the resulting curves are shown in 
Figures~\ref{fig:experiments:Dirichlet:Diaphragm:boundary_data} and~\ref{fig:experiments:Diaphragm:Robin:boundary_data}.

\section{\victor{Benchmark I: Deformation of the diaphragm using the smoothed Dirichlet boundary conditions}}
\label{sec:experiments:dirichlet}
In this section we solve the discretized linear elasticity equations \eqref{eq:method:DandF}. 
We are interested in whether the problem with smoothed Dirichlet 
boundary condition 
leads to a high-order convergence. A point of interest is also whether the resulting deformation is 
physiologically sensible, and whether the Von Mises stress is distributed as expected.

The boundary condition that we use is purely Dirichlet, which means
that we set the Robin coefficients in \eqref{eq:model:PDE} 
to $\kappa_0(y) = 1$ and $\kappa_1(y)=0$. Then the only boundary data functions present in the system 
\eqref{eq:method:DandF} are $g_1(y)$ and $g_2(y)$, which correspond 
to the imposition of displacements in the horizontal and vertical direction, respectively.

\subsection{The imposition of boundary displacements}
\begin{figure}[htb!]
    \centering
    \begin{tabular}{@{}cc@{}}
        \hspace{1.3cm}\textbf{Displacement $\mathbf{\tilde u_1}$} & \hspace{1.3cm}\textbf{Displacement $\mathbf{\tilde u_2}$} \\
        \includegraphics[height=0.22\linewidth]{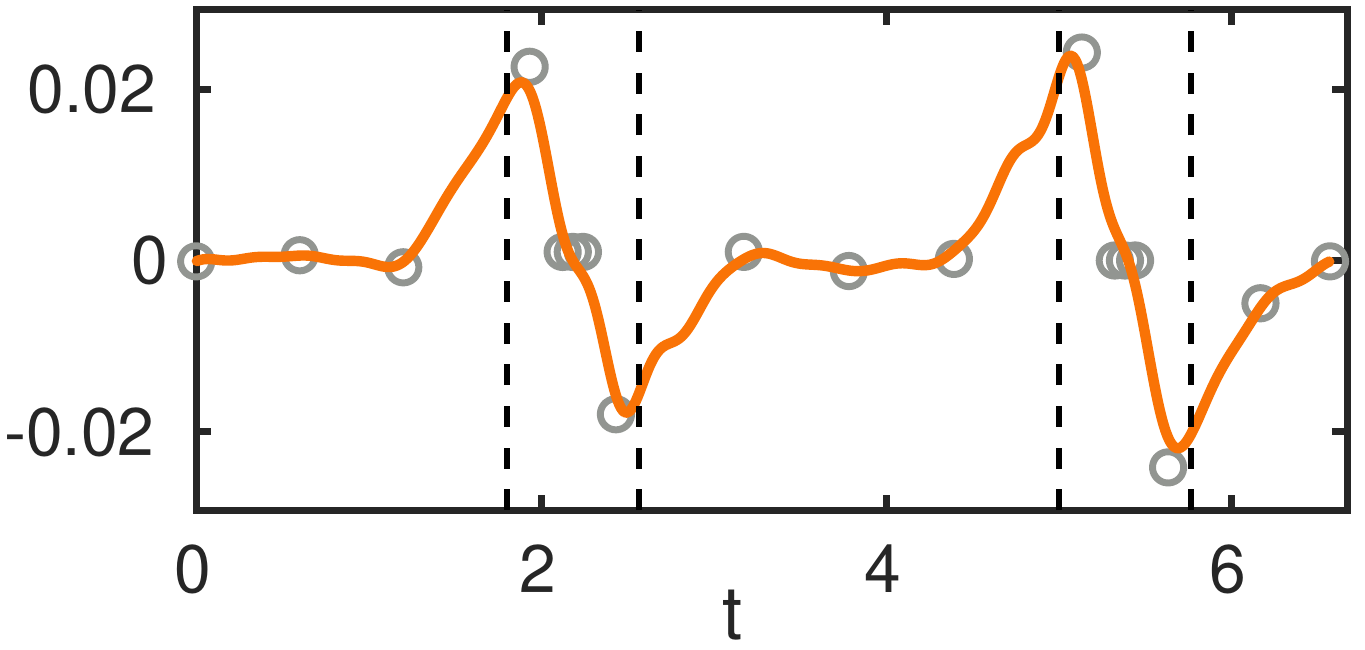} & 
        \raisebox{0.005\linewidth}{\includegraphics[height=0.22\linewidth]{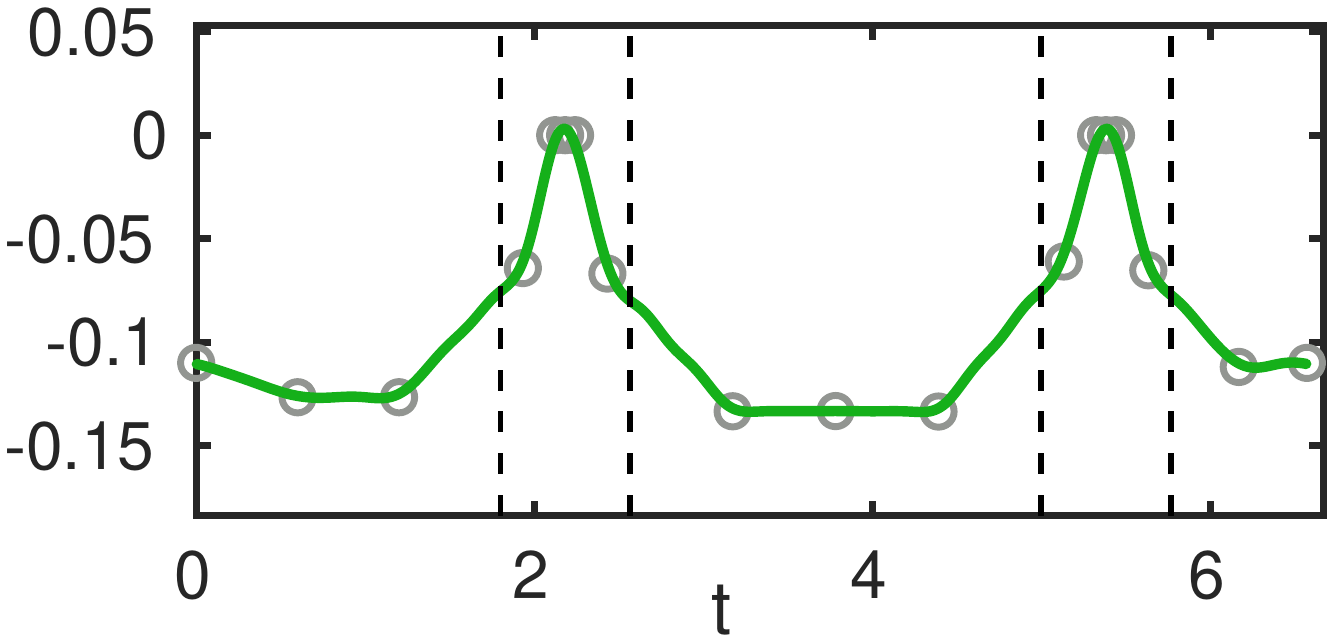}}
    \end{tabular}
    \caption{Benchmark I: The displacement in the horizontal direction (left) and the vertical direction (right) as a 
    function of the boundary parametrization $t$ (see Figure~\ref{fig:boundary}, right image, for an illustration of $t$ in relation to the boundary). The markers show the initially placed data points, which are then linearly interpolated into $\tilde{M}_g=80$ data points. For the approximation, $N_g=120$ node points and stencil size $n=28$ were used. 
    The dashed lines show the location of the end points of the diaphragm (regions 2 and 8 in Figure~\ref{fig:smoothing:bcPosition}).}
    \label{fig:experiments:Dirichlet:Diaphragm:boundary_data}
\end{figure}

\begin{figure}[htb!]
    \centering
    \includegraphics[width=0.6\linewidth]{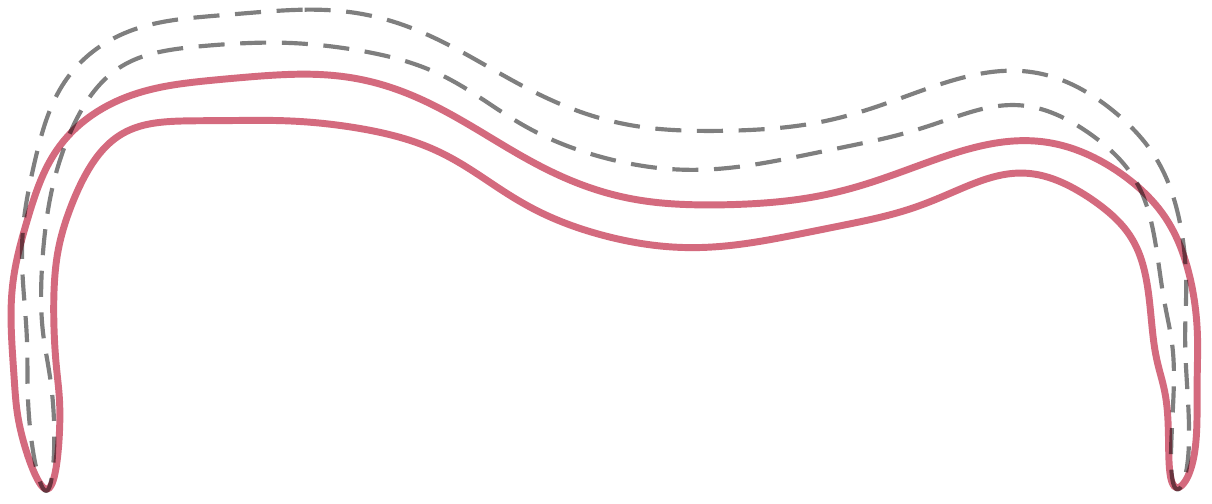}
    \caption{Benchmark I: The diaphragm in its non-deformed state (dashed line) and after displacement of the boundary (solid line).}
    \label{fig:experiments:Dirichlet:Diaphragm:deformed}
\end{figure}

%
%
The displacements have been synthesized to reproduce the physiological behavior described in Section 
\ref{section:diabehavior}. Particularly, the translation of the horizontal part of the diaphragm (region 5) has been defined as a constant vertical displacement 
in the downward direction and the thickening of the appositional zone (regions 3 and 5) is also prescribed as a constant. 
As described in Section~\ref{sec:smooth}, we place a few data points based on this information, and then the rest of the data is 
generated by linear interpolation. From a physiological perspective all displacements should be smooth. 
We generate a smooth function by solving the constrained least squares problem~\eqref{eq:boundary:system}. 
%
%
The results are shown in Figure~\ref{fig:experiments:Dirichlet:Diaphragm:boundary_data}.
Imposing smooth boundary data makes it possible to obtain high-order convergence and provides a solution that is physically relevant. 
In Figure~\ref{fig:experiments:Dirichlet:Diaphragm:deformed} we display the boundary of the diaphragm before and after 
application of the displacements from Figure~\ref{fig:experiments:Dirichlet:Diaphragm:boundary_data}.
\subsection{Solution of Benchmark I}
\begin{figure}[h!]
    \centering    
    \begin{tabular}{c}
    \textbf{Displacement $\tilde u_1$} \\
    \includegraphics[width=0.65\linewidth]{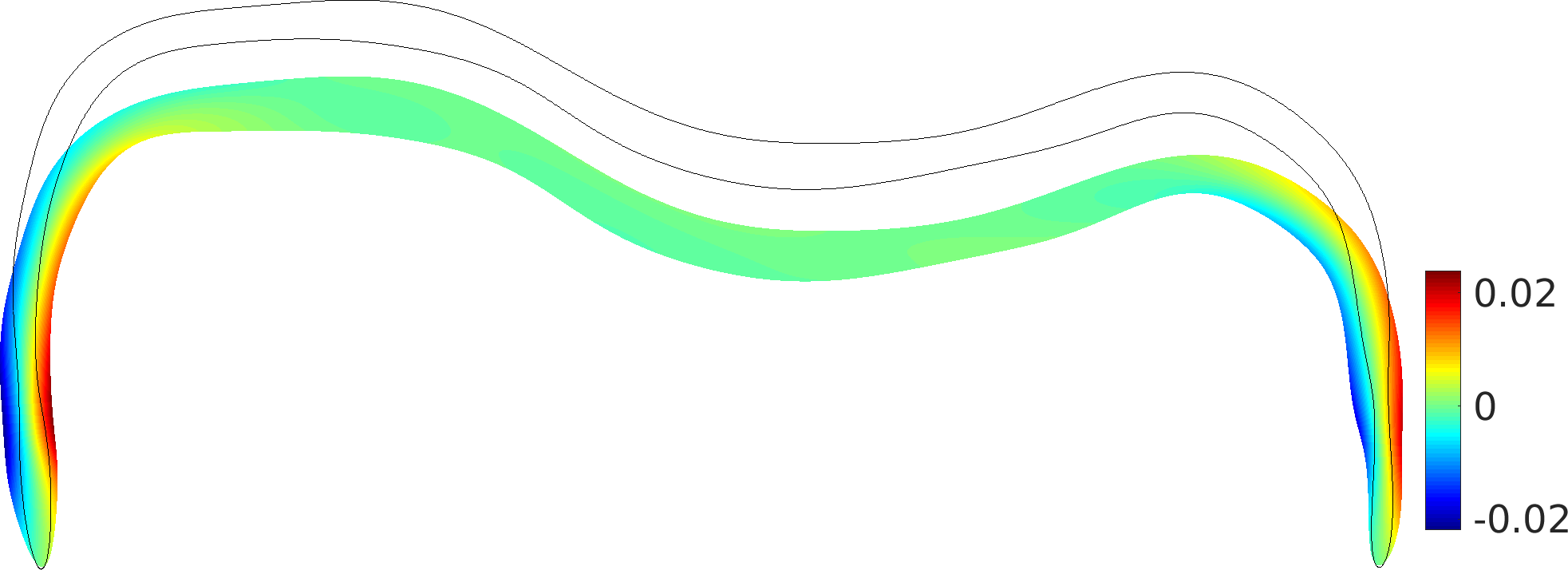} \\
    \textbf{Displacement $\tilde u_2$} \\
    \includegraphics[width=0.65\linewidth]{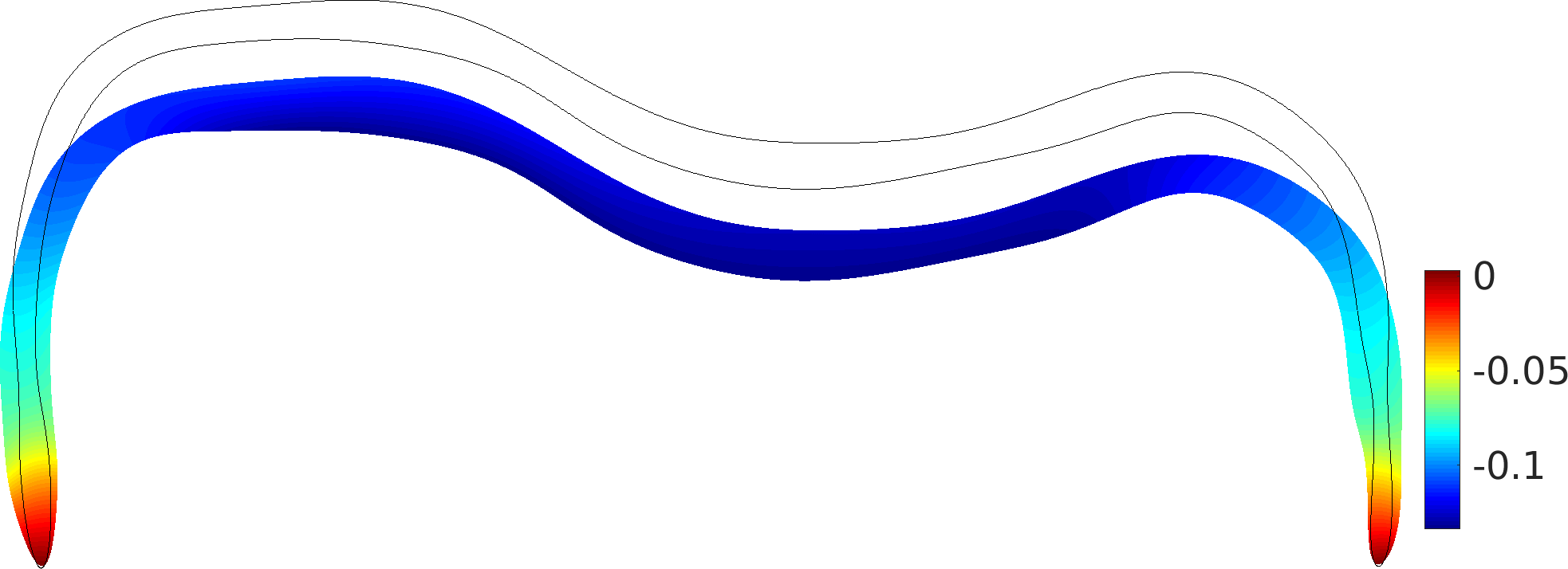} \\
    \textbf{Von Mises stress} \\
    \includegraphics[width=0.65\linewidth]{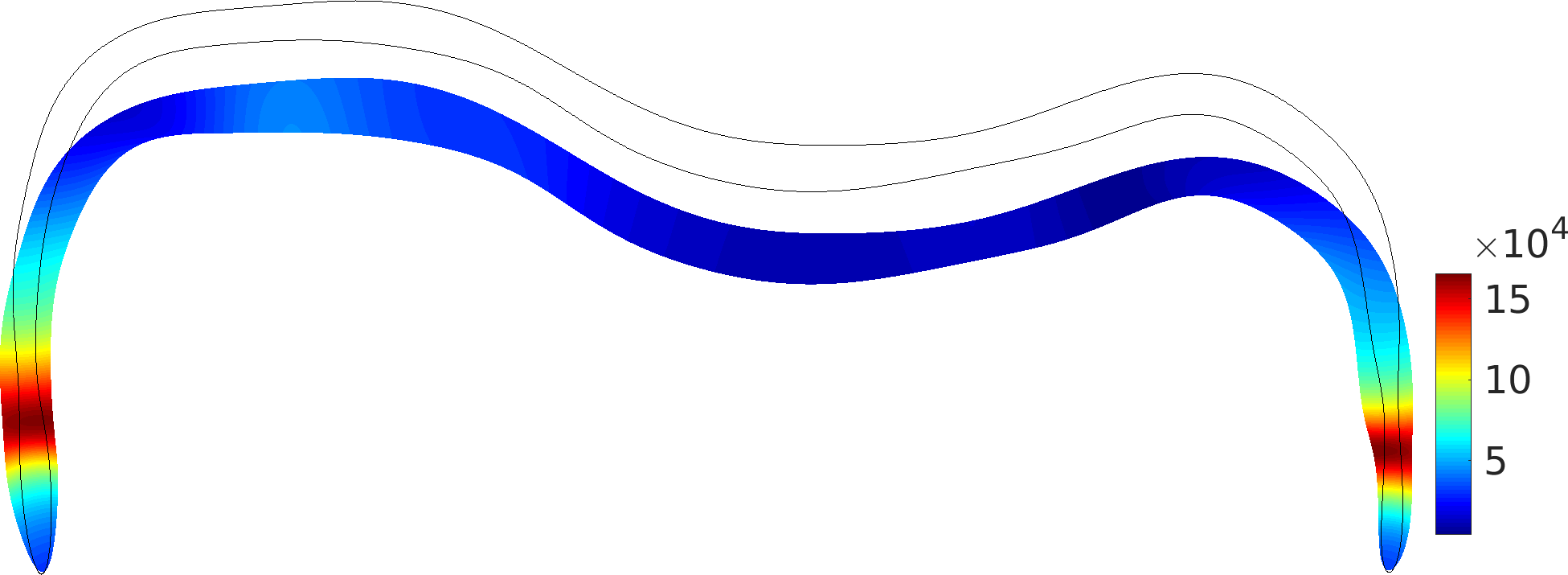}
    \end{tabular}
    \caption{Benchmark I:  The computed displacements and the corresponding von Mises stress over the diaphragm. This solution was obtained using the unfitted RBF-FD discretization with internodal distance $h=0.004$, 
    oversampling parameter $q=5$, and an appended polynomial basis of degree $p=5$. \victor{Due to the scaling applied to the geometry (see Section~\ref{sec:smooth}), the displayed results for displacement and stress 
    should be multiplied with $s_{\Omega}^{-1} = 0.15692$ m to recover the results corresponding to the unscaled problem.}}
    \label{fig:experiments:Dirichlet:Diaphragm:solution}
\end{figure}

The solution is given in Figure~\ref{fig:experiments:Dirichlet:Diaphragm:solution}, 
where we display the spatial distribution of the displacements and the Von Mises stress. 
Looking at the $\tilde u_1$ distribution in the figure, we can observe thickening of the diaphragm 
in regions 2, 3, 7 and 8 according to the labels from Figure~\ref{fig:smoothing:bcPosition}. 
Next, the distribution of $\tilde u_2$ indicates that translation is largest in regions 4, 5 and 6. 
The Von Mises stress is largest at the interfaces between regions 2, 3 and regions 7, 8. \elisabeth{This makes sense, since the change in the thickness is largest in these regions}. 
We conclude that the behavior roughly follows the 
physiological cues described in Section~\ref{section:diabehavior}.

\subsection{Convergence under node refinement}
While the solution is roughly what we expect from a physiological perspective, we are yet to understand whether 
the simulation gives a correct answer from a numerical perspective. 
We investigate the convergence of the numerical solution under node refinement using several different 
polynomial degrees in the stencil-based approximation.
Since we do not know what the true solution is, we measure convergence of the numerical solution $\tilde u(Y)$ towards a numerical reference 
solution $\tilde u_*(Y_*)$, where \elisabeth{the node set is highly refined}. We choose two numerical references: 
(i) computed using the unfitted RBF-FD-LS method with internodal distance $h=0.002$, leading to $N=43\,841$ and polynomial degree $p=5$,  
(ii) computed using the Galerkin finite element method with 
linear elements and $236\,414$ degrees of freedom (corresponding to $h = 0.00085$), using the GetDP solver \cite{getdp-siam2008}. 
Every numerical solution $\tilde u$ that we obtain is interpolated, \victor{consistent with the approximation order,} to the point set 
of the reference solution $Y_*$, where we then compute the approximation error:
\begin{equation}
    \label{eq:dirichlet:error}
\| e \|_{\ell_2} = \frac{\|\tilde u(Y_*) - \tilde u_*(Y_*)\|_2}{\|\tilde u_*(Y_*)\|_2}.
\end{equation}
\begin{figure}[h!]
    \centering
    \begin{tabular}{ccc}
        \hspace{0.9cm}\textbf{Displacement } $\mathbf{\tilde u_1}$ & \hspace{0.9cm}\textbf{Displacement } $\mathbf{\tilde u_2}$ & \hspace{0.9cm}\textbf{Von Mises} \\
        \includegraphics[width=0.315\linewidth]{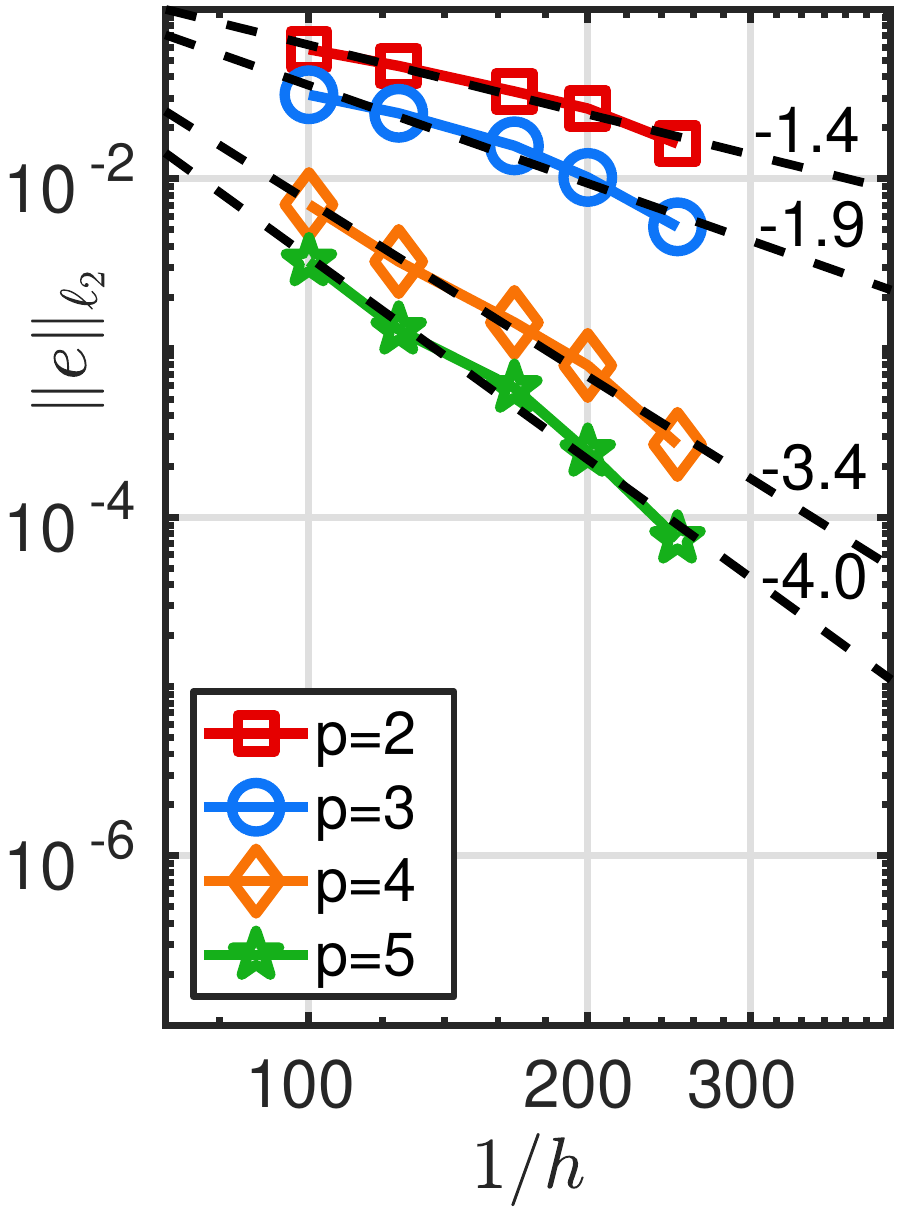} &  
        \includegraphics[width=0.315\linewidth]{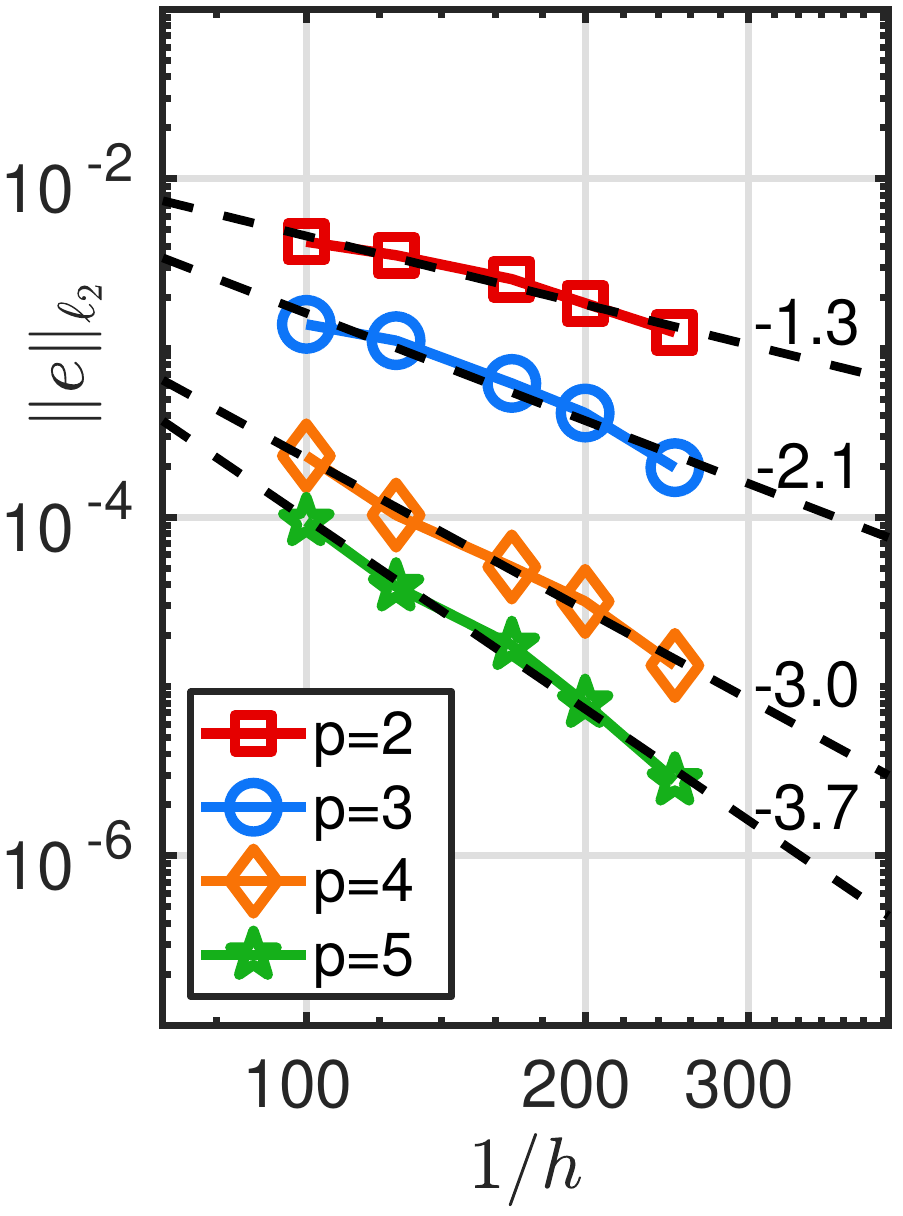} &  
        \includegraphics[width=0.315\linewidth]{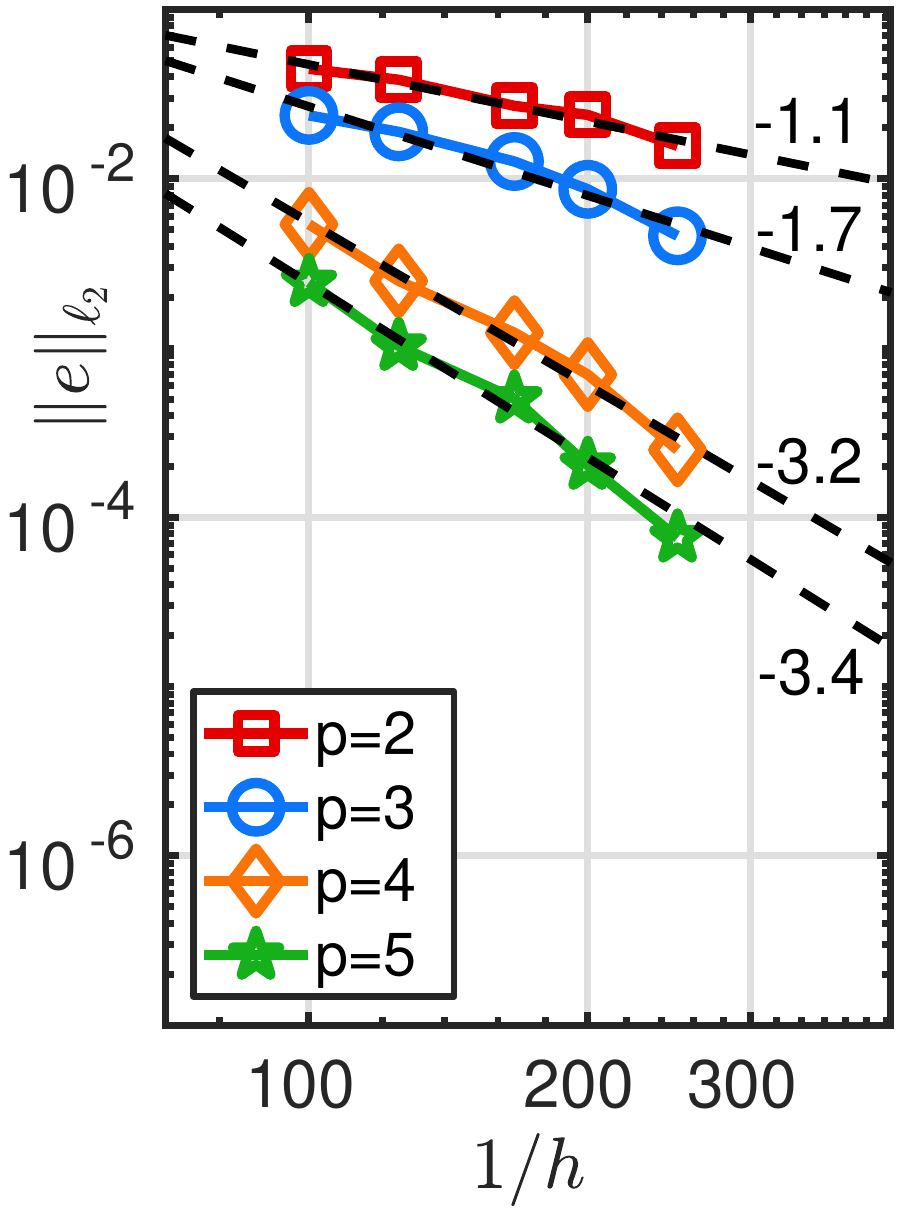}  
    \end{tabular}        
    \caption{Benchmark I: Convergence of the displacements $\tilde u_1$ and $\tilde u_2$ and the Von Mises stress for different polynomial degrees $p$, 
    against a highly resolved numerical solution computed using the unfitted RBF-FD-LS method, \victor{with $h=0.002$ and $p=5$}.}
    \label{fig:experiments:Dirichlet:Diaphragm:convergence_selfRef}
\end{figure}
Table~\ref{fig:experiments:Dirichlet:Diaphragm:h_to_N_table} shows the relation between the 
internodal distance $h$, computed according to \eqref{eq:method:h_and_hy}, and the number of degrees of freedom $N$ for the considered problem sizes. 
\begin{table}[h!]
    \centering
    \caption{Benchmark I and II: The relation between the inverse internodal distance $1/h$ and (i) the internodal distance $h$ and (ii) the number of interpolation points $N$ used for discretizing the PDE problem \eqref{eq:model:PDE}.}
    \begin{tabular}{|c|c|c|c|c|c|c|c|c|}
      \hline
        $\mathbf{1/h}$ & 25 & 50 & 100 & 125 & 166.67 & 200 & 250 & 500 \\ \hline
        $\mathbf{h}$ &  0.04 & 0.02 & 0.01 & 0.008 & 0.006 & 0.005 & 0.004 & 0.002 \\ 
        $\mathbf{N}$ & 337 & 853 & 2500 & 3617 & 5904  & 8126 & 12129 & 43841 \\
        \hline
    \end{tabular}        

    \label{fig:experiments:Dirichlet:Diaphragm:h_to_N_table}
\end{table}

The results when the unfitted RBF-FD method is used as a reference are displayed in Figure~\ref{fig:experiments:Dirichlet:Diaphragm:convergence_selfRef}. 
We observe that the error is small for all polynomial degrees. \elisabeth{The convergence rate increases as $p$ is increased.} 
For all $p$ the convergence rates of the displacements are close to $p-1$. 
The convergence rate of the Von Mises stress is larger than $p-2$ for every $p$. 
Stress is computed using the first derivatives of the displacements, which (in theory) lowers the 
convergence rate with one \elisabeth{order}.

\begin{figure}[h!]
    \centering
    \begin{tabular}{ccc}
        \hspace{0.9cm}\textbf{Displacement }$\mathbf{\tilde u_1}$ & \hspace{0.9cm}\textbf{Displacement }$\mathbf{\tilde u_2}$ & \hspace{0.9cm}\textbf{Von Mises} \\
        \includegraphics[width=0.315\linewidth]{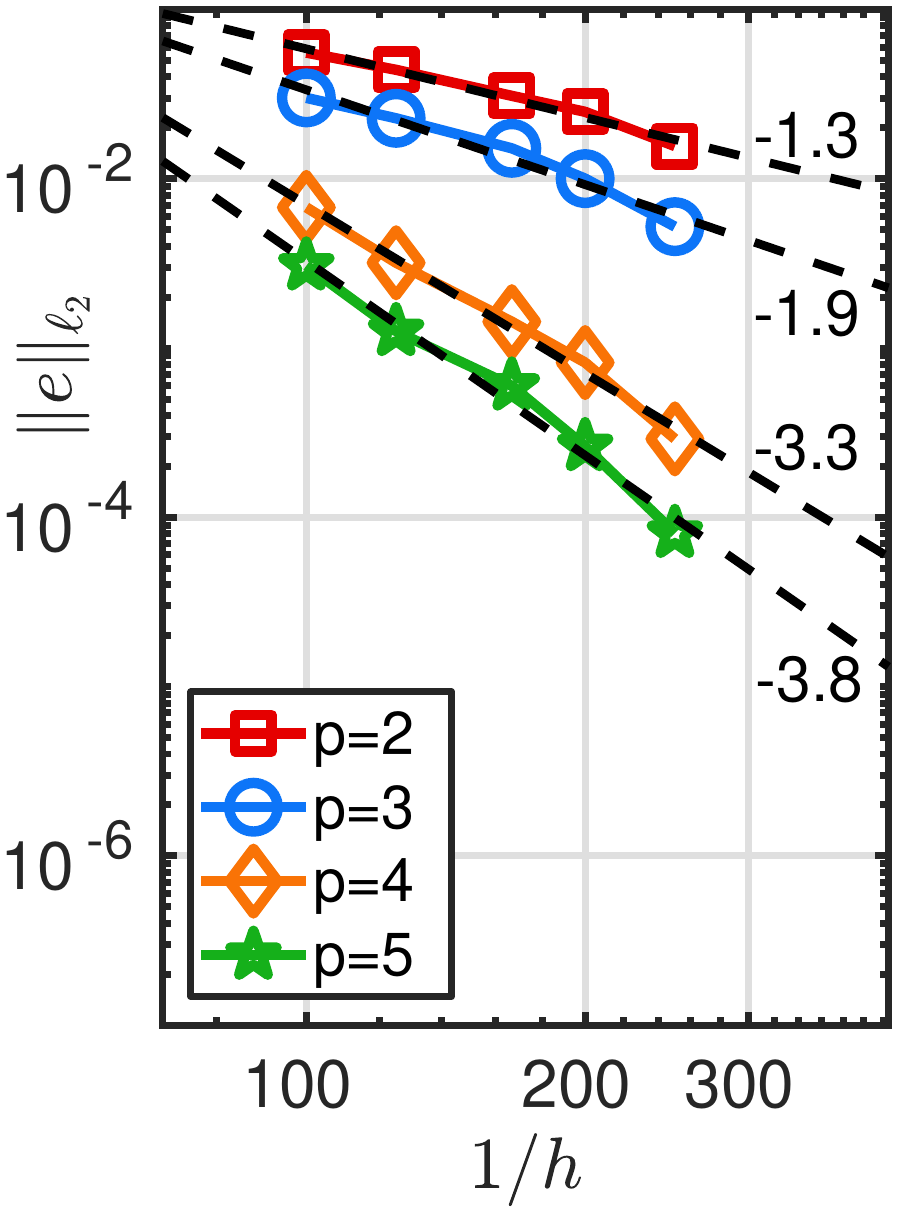} &  
        \includegraphics[width=0.315\linewidth]{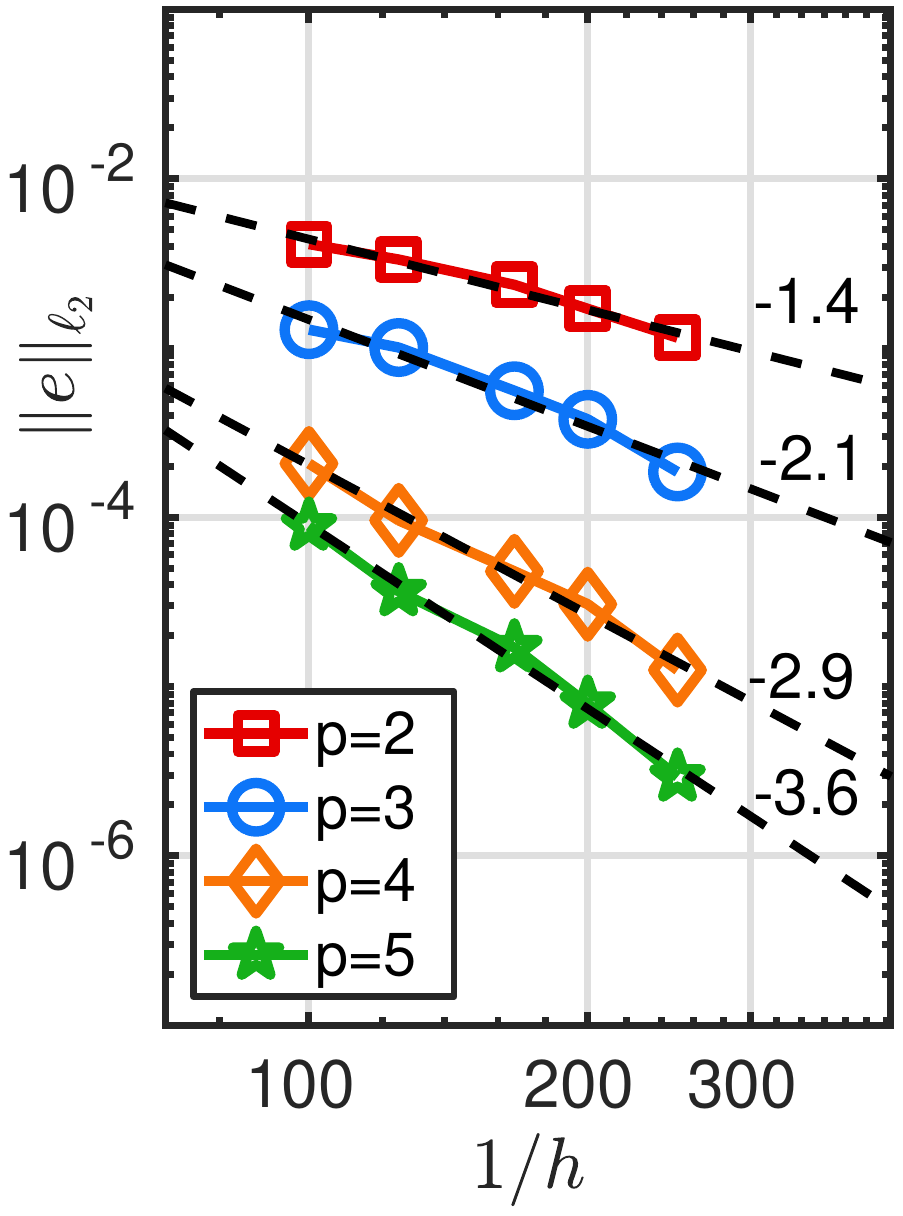} &  
        \includegraphics[width=0.315\linewidth]{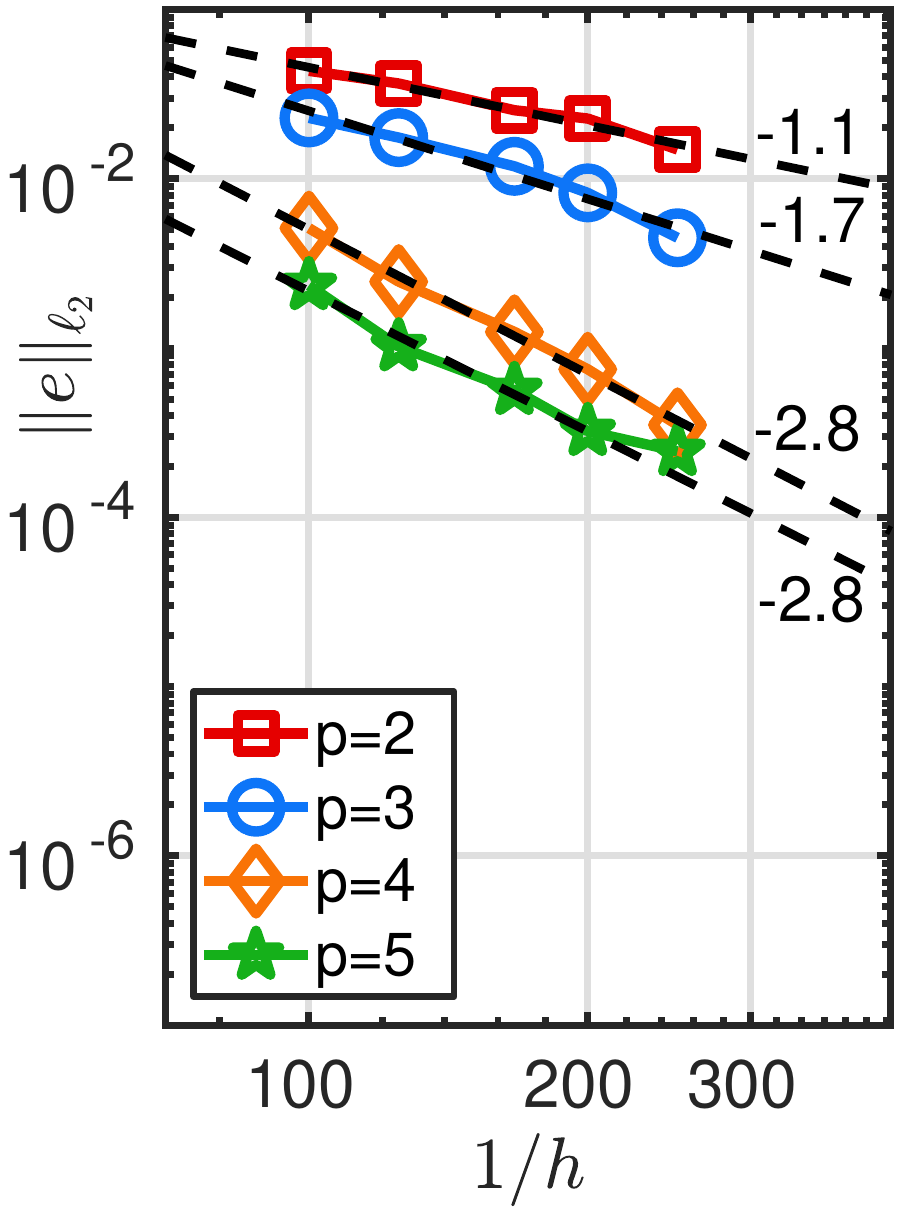}  
    \end{tabular}    
    \caption{Benchmark I: Convergence of the displacements $\tilde u_1$ and $\tilde u_2$ and the Von Mises stress for different polynomial degrees $p$, 
    against a highly resolved numerical solution computed using FEM \victor{(the GetDP solver \cite{getdp-siam2008}) with linear elements and $236\,414$ degrees of freedom (corresponding to $h = 0.00085$).}}
    \label{fig:experiments:Dirichlet:Diaphragm:convergence_femRef}
\end{figure}

The results when the finite element method is used as a reference are given in Figure~\ref{fig:experiments:Dirichlet:Diaphragm:convergence_femRef}. 
The convergence plots for the displacement are very similar to the results in Figure~\ref{fig:experiments:Dirichlet:Diaphragm:convergence_selfRef}.
The convergence rate of the stress for $p=5$ is \elisabeth{lower} compared with the self-reference test provided in Figure~\ref{fig:experiments:Dirichlet:Diaphragm:convergence_selfRef}. 
Furthermore, when $p=5$, the convergence seems to be stalling at the last point of observation. 
Our speculative reasoning is that the derivatives in the finite element space 
do not approximate the stresses well enough there. 
This implies that the FEM mesh is too coarse to match the accuracy of the unfitted RBF-FD method for derivative approximation, 
when the polynomial degree is as high as $p=5$, for the range of $h$ used here.

\begin{figure}[h!]
    \centering
    \begin{tabular}{cc}
    \textbf{Displacement $\tilde u_1$} & \textbf{Displacement $\tilde u_2$} \\
    \includegraphics[width=0.49\linewidth]{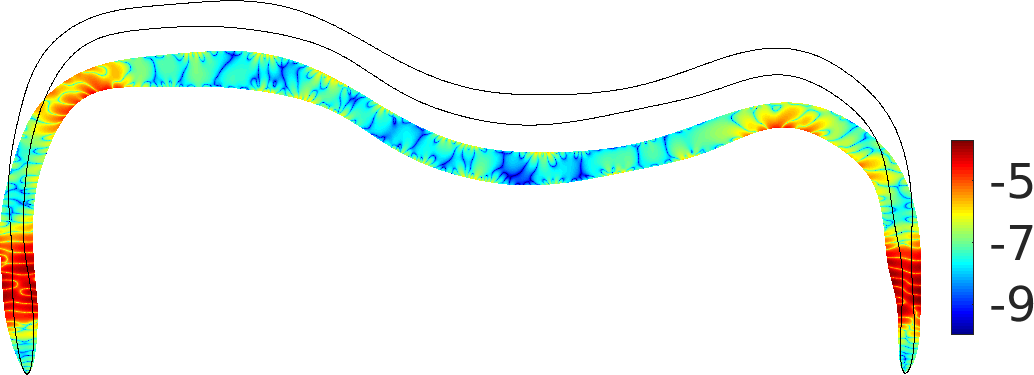} &
    \includegraphics[width=0.49\linewidth]{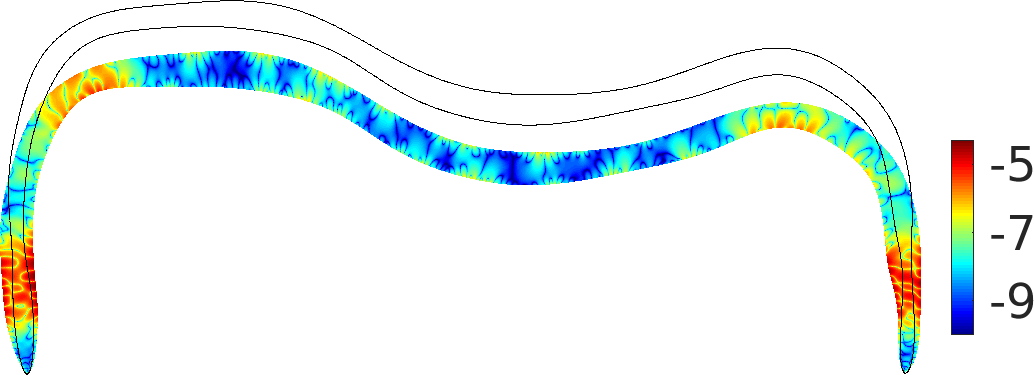} \\
    \multicolumn{2}{c}{\textbf{Von Mises stress}} \\
    \multicolumn{2}{c}{\includegraphics[width=0.6\linewidth]{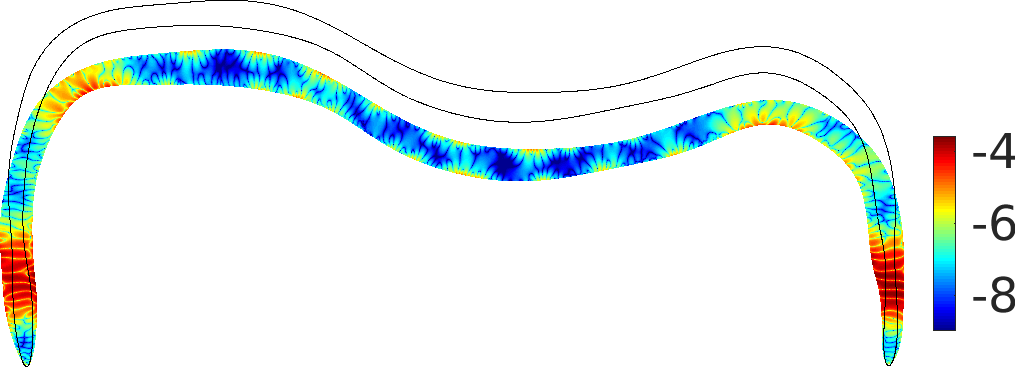}}
    \end{tabular}
    \caption{Benchmark I: Error distribution in logarithmic scale for $h=0.006$ and $p=5$ when a fine unfitted RBF-FD-LS solution ($h=0.002$, $p=5$) 
    is taken as reference for computing the error.}
    \label{fig:experiments:Dirichlet:Diaphragm:error_spatial_selfRef}
\end{figure}

\begin{figure}[h!]
    \centering
    \begin{tabular}{cc}
    \textbf{Displacement $\tilde u_1$} & \textbf{Displacement $\tilde u_2$} \\
    \includegraphics[width=0.49\linewidth]{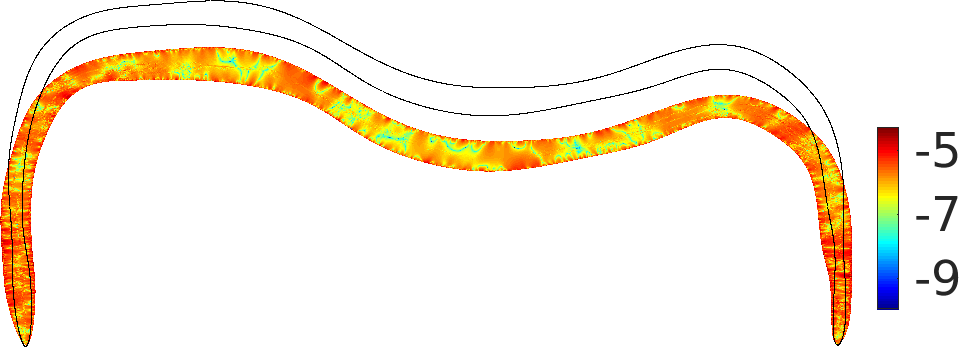} &
    \includegraphics[width=0.49\linewidth]{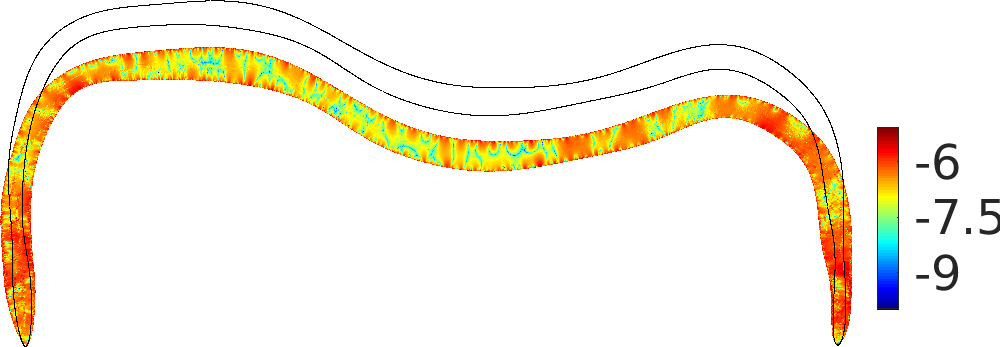} \\
    \multicolumn{2}{c}{\textbf{Von Mises stress}} \\
    \multicolumn{2}{c}{\includegraphics[width=0.6\linewidth]{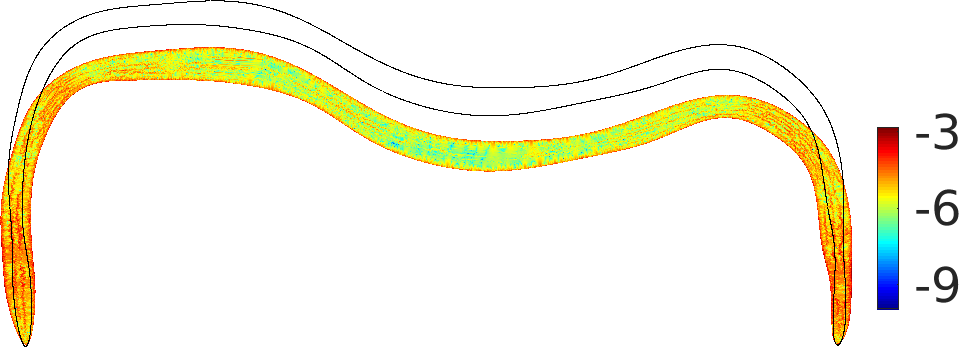}}
    \end{tabular}
    \caption{Benchmark I: Error distribution in logarithmic scale for $h=0.006$ and $p=5$ when a fine unfitted RBF-FD-LS solution ($h=0.002$, $p=5$) 
    is taken as reference for computing the error.}
    \label{fig:experiments:Dirichlet:Diaphragm:error_spatial_femRef}
\end{figure}

The spatial distribution of the error when using the self-reference solution and the finite element reference solution 
is shown in Figure 
\ref{fig:experiments:Dirichlet:Diaphragm:error_spatial_selfRef} and Figure~\ref{fig:experiments:Dirichlet:Diaphragm:error_spatial_femRef}, 
respectively. The numerical solution was computed using $h=0.004$, $q=5$, $p=5$.
\elisabeth{For all solution fields and both references we can observe that the error is larger in the regions with larger Von Mises stress. In addition, the error also tends to be larger in the regions where 
the boundary curve is concave. This implies that our problem could benefit from adaptive 
node refinement, which we are planning to use in our future work.}

\section{\victor{Benchmark II: Deformation of the diaphragm using the smoothed Robin boundary conditions}}
\label{sec:experiments:Robin}
This benchmark includes a more difficult problem compared with Benchmark I from Section 
\ref{sec:experiments:dirichlet}. As discussed in Section \ref{section:diabehavior}, we can measure 
the displacements over certain parts of the boundary, i.e., where the diaphragm is fixed near the spine, or where it moves together 
with the sternum and some ribs. In the regions where the displacements are not known, 
we may instead have information about the thoracic or abdominal pressure. A pressure condition is a special case of a traction boundary condition. 
A straightforward way to handle these boundary conditions would be to impose the Dirichlet condition (the known displacements) 
and the traction condition (the known traction values) in disjoint regions. This implies a discontinuous imposition of boundary conditions, 
for which: (i) we can not obtain high-order convergence according to our preliminary tests, (ii) 
the resulting deformation might be discontinuous or have large local derivatives, which would not reflect 
the physicological behavior of the diaphragm. For this reason we introduce a smooth blending of the Dirichlet and the traction 
boundary conditions in Robin form~\eqref{eq:model:secondNewton_bcs_Robin}.

\begin{figure}[h!]
    \centering
    \includegraphics[width=0.65\linewidth]{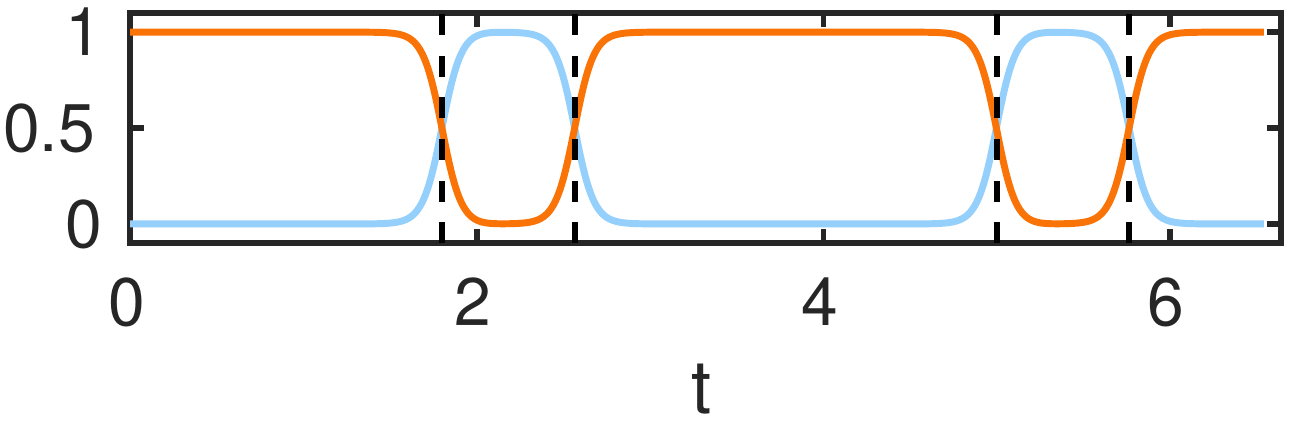}

    \caption{The Robin boundary coefficients displayed over the diaphragm show the impositions of the traction 
    and the Dirichlet parts of the boundary condition. The orange line corresponds to the traction coefficient $\kappa_1$, 
    and the blue line corresponds to the Dirichlet coefficient $\kappa_0$. The dashed lines show the location of the end 
    points of the diaphragm (regions 2 and 8 in Figure~\ref{fig:smoothing:bcPosition}). The parameter $t$ corresponding 
    to the boundary of the diaphragm is illustrated in Figure \ref{fig:boundary}.}

    \label{fig:experiments:Diaphragm:Robin:coefficients}
\end{figure}

This setting implies that we solve the system of equations \eqref{eq:method:DandF},
where the supports of the Robin coefficients $\kappa_0(y)$ and $\kappa_1(y)$ are chosen to overlap slightly in the regions where the type of boundary data changes.
In Figure \ref{fig:experiments:Diaphragm:Robin:coefficients} we display the Robin coefficients which we use in Benchmark II. 
The displayed coefficients are a function of the boundary parameter $t$. 
The coefficients were computed using a sum of the sigmoid functions $\frac{1}{1+e^{\varepsilon (t-d_i)}}$, where $\varepsilon = \pm 20$, 
$t$ is the boundary parameter and $d_i$, $i=1,2,..$ are the transition points 
(marked by dashed lines in Figure \ref{fig:experiments:Diaphragm:Robin:coefficients}). The 
sign of $\varepsilon$ depends on whether the coefficient is increasing (positive sign) or decreasing (negative sign).

\begin{figure}[h!]
    \centering
    \begin{tabular}{@{}c@{}cc@{}}
      &\hspace*{8mm}\textbf{Dirichlet data} & \hspace*{4mm}\textbf{Traction data} \\
      \rotatebox{90}{\hspace*{4mm}\textbf{Displacements }$\mathbf{\tilde u_1}$}  
        &\includegraphics[width=0.48\linewidth]{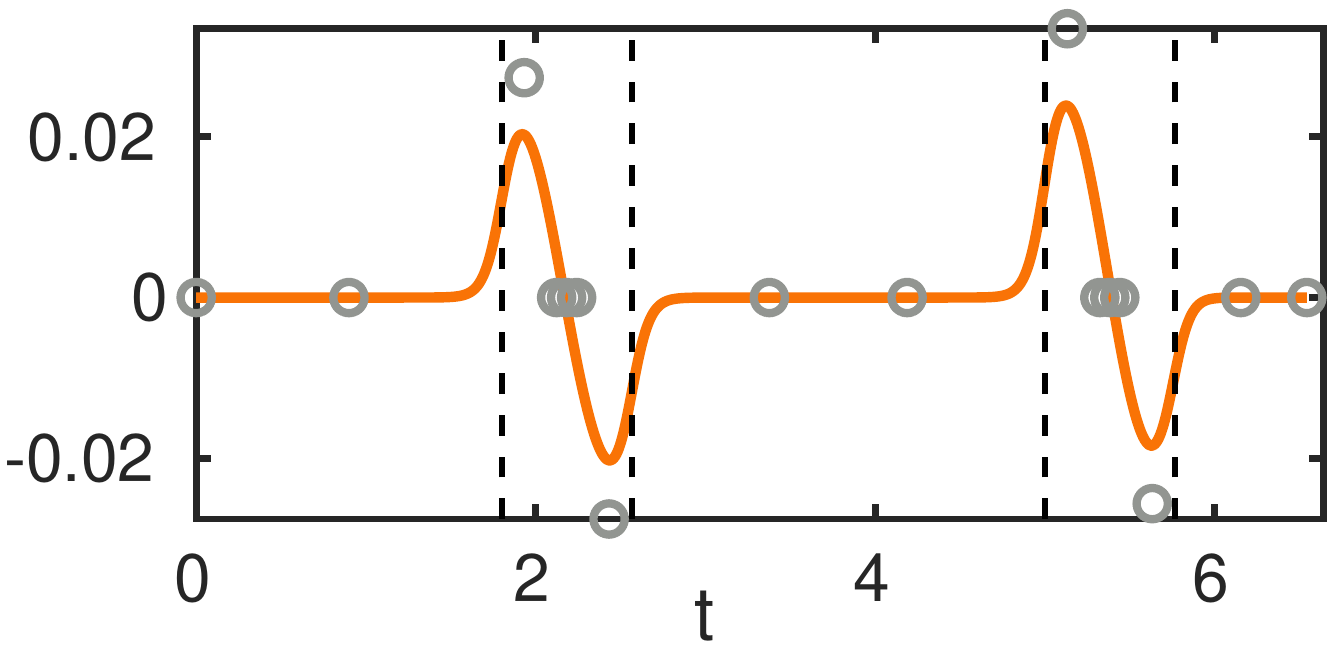} & 
        \raisebox{0.005\linewidth}{\includegraphics[width=0.445\linewidth]{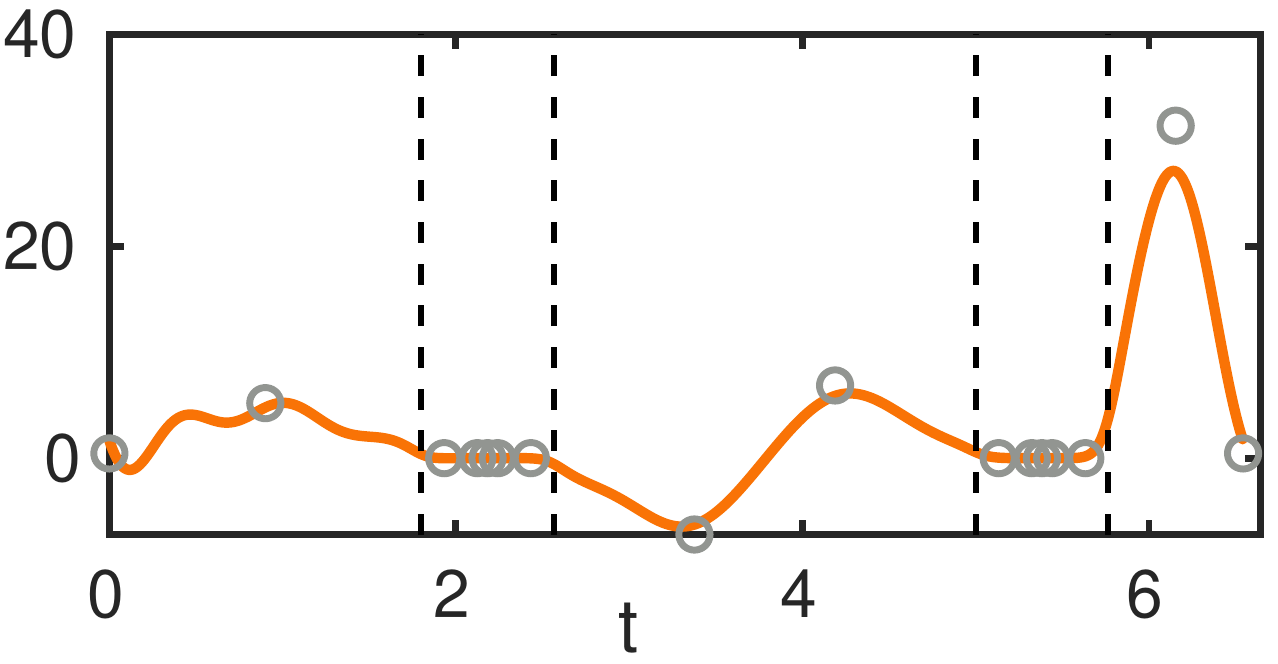}} \\
        \rotatebox{90}{\hspace*{4mm}\textbf{Displacements }$\mathbf{\tilde u_2}$}
        &\hspace{0.04\linewidth}\includegraphics[width=0.45\linewidth]{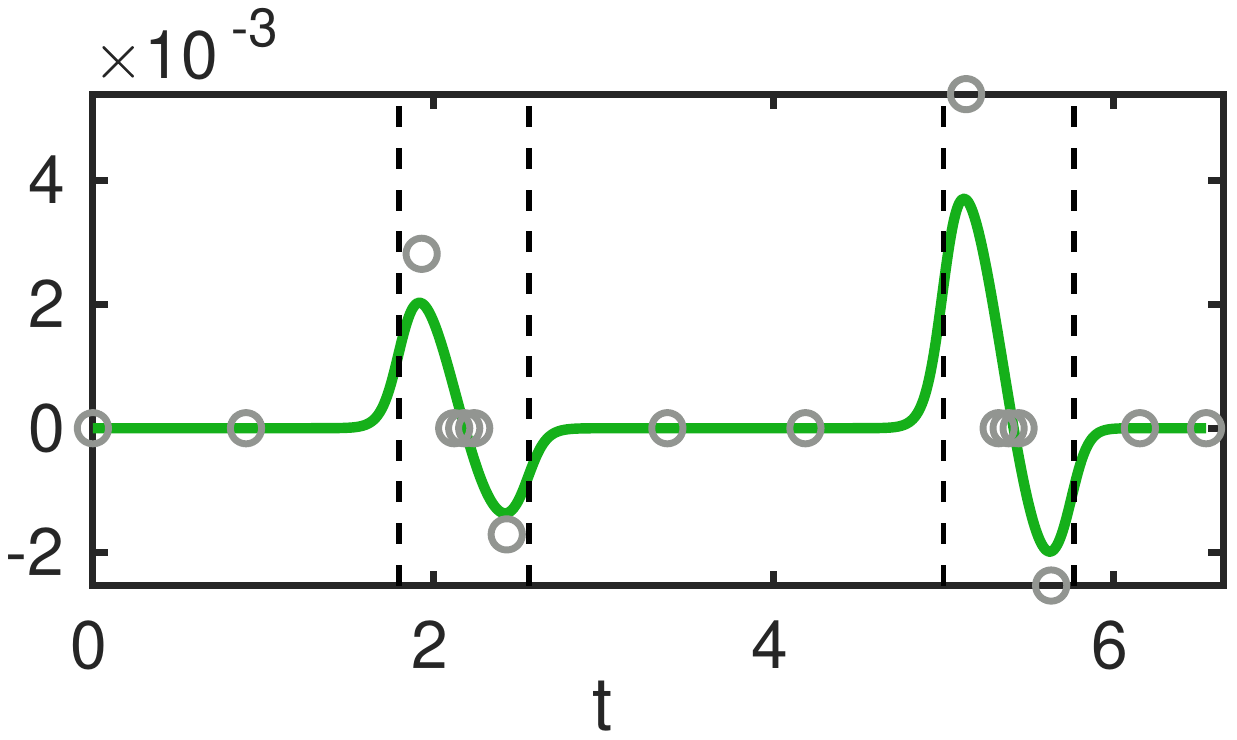} &
        \raisebox{0.01\linewidth}{\includegraphics[width=0.465\linewidth]{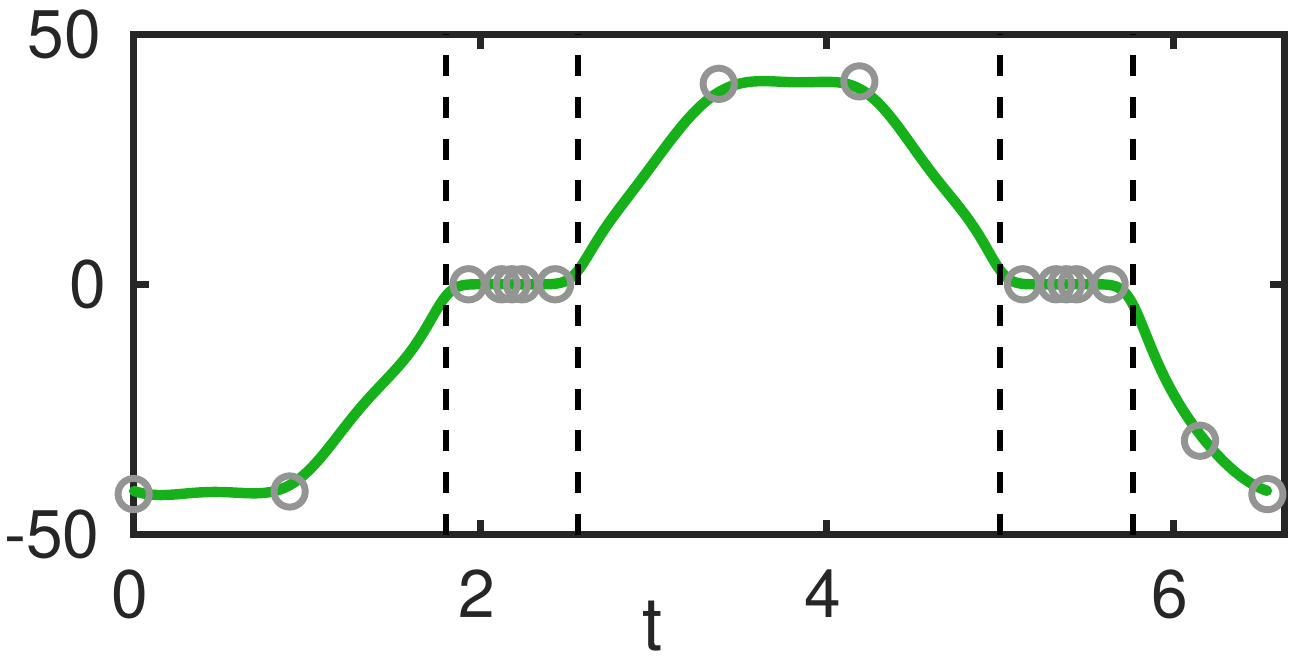}}
    \end{tabular}
    \caption{Benchmark II: The smoothed boundary data over the diaphragm. For this benchmark we have both Dirichlet data (left column) 
    and traction data (right column), each with a horizontal component (first row) and a vertical component (second row). 
    The markers show the initially placed data points, which are then linearly interpolated into $\tilde{M}_g=80$ data points. 
    For the approximation, $N_g=120$ node points and stencil size $n=28$ were used. 
    \victor{The dashed lines show the location of the end points of the diaphragm (regions 2 and 8 in Figure~\ref{fig:smoothing:bcPosition}). 
    The parameter $t$ corresponding to the boundary of the diaphragm 
    is illustrated in Figure \ref{fig:boundary} (right image).}}
    \label{fig:experiments:Diaphragm:Robin:boundary_data}
\end{figure}

The smoothed boundary displacement values $g_1$, $g_2$ and the boundary traction values $h_1$, $h_2$ are given in Figure \ref{fig:experiments:Diaphragm:Robin:boundary_data}. 
Through the functions $g_1$ and $g_2$ we impose thickening in the regions 2, 3, 7 and 8 based on the distribution in Figure \ref{fig:smoothing:bcPosition}. 
Additionally, we also impose a gentle translation in the negative vertical direction in the regions 7 and 8. 
Through the functions $h_1$ and $h_2$ we impose the translation of zone 5 in the negative vertical direction.
\subsection{Solution of Benchmark II}
\begin{figure}[h!]
    \centering
    
    \begin{tabular}{c}
    \textbf{Displacement $\mathbf{\tilde u_1}$} \\
    \includegraphics[width=0.65\linewidth]{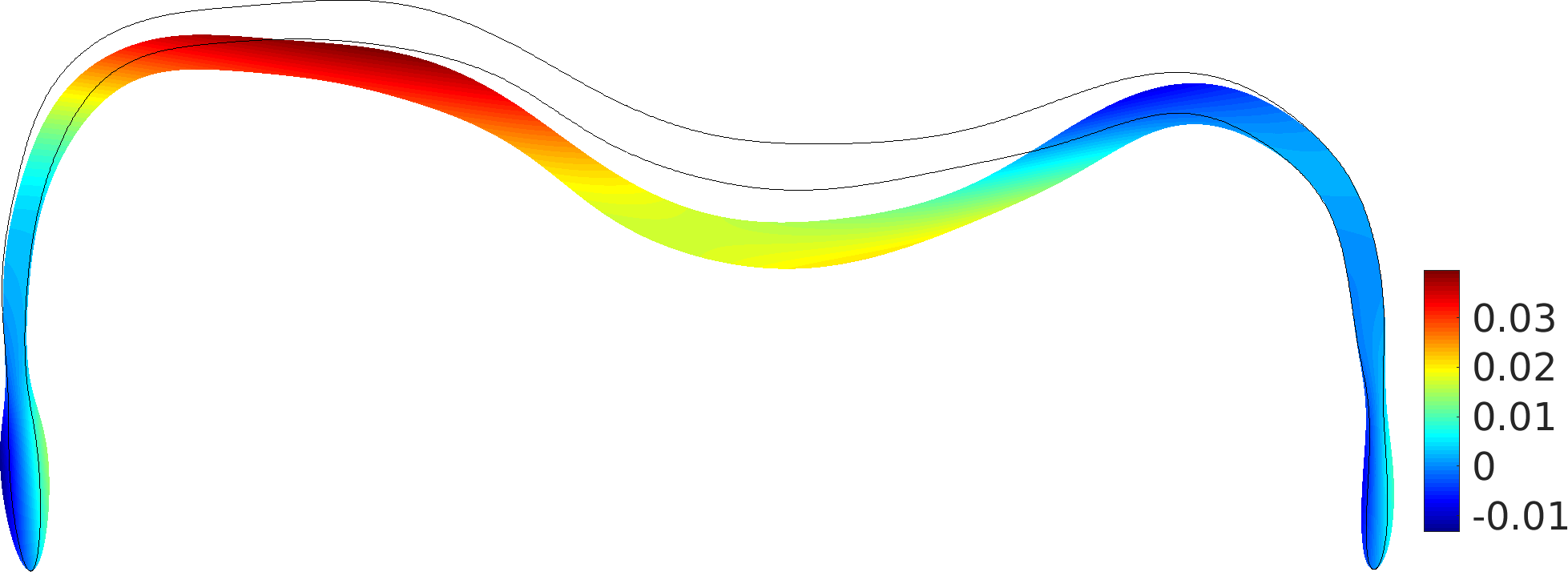} \\
    \textbf{Displacement $\mathbf{\tilde u_2}$} \\
    \includegraphics[width=0.65\linewidth]{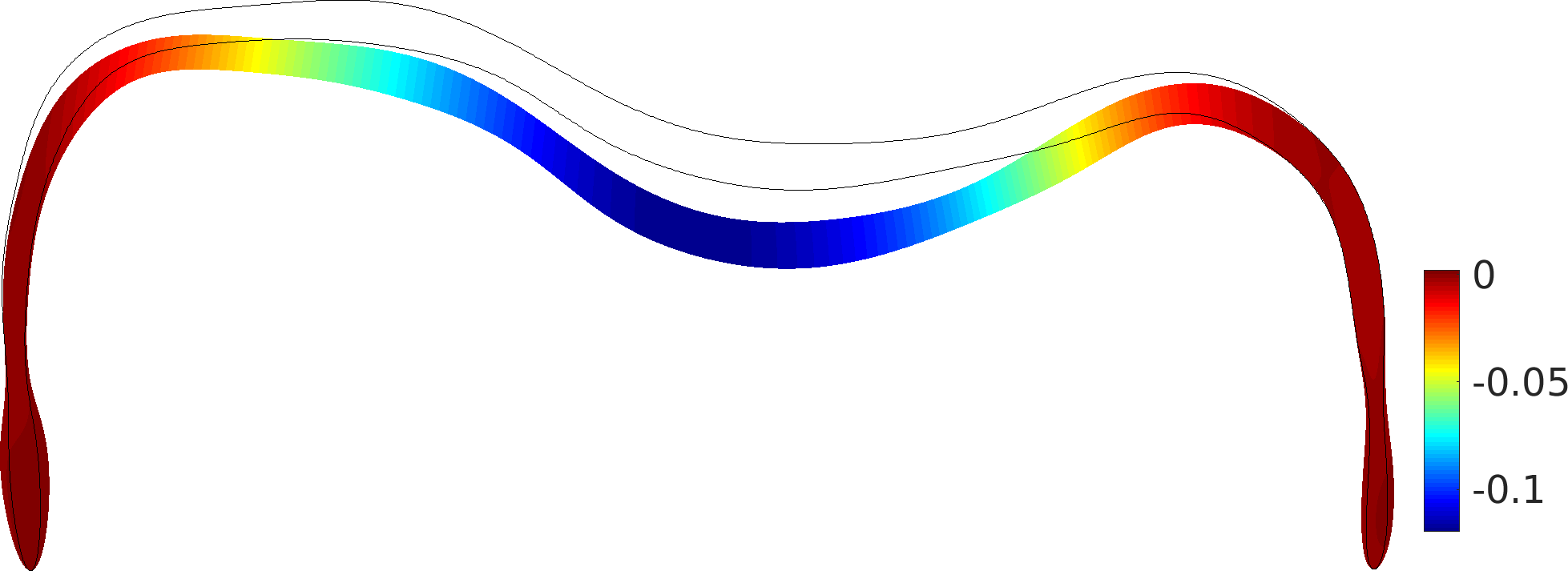} \\
    \textbf{Von Mises stress} \\
    \includegraphics[width=0.65\linewidth]{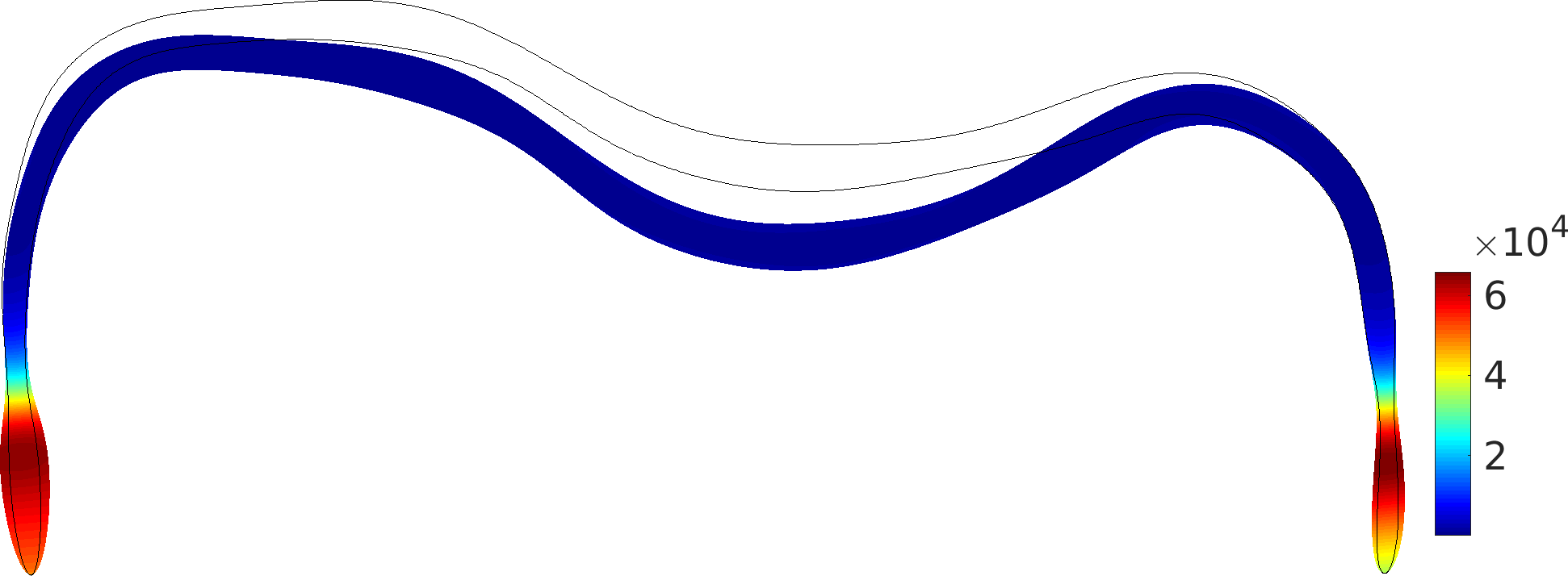}
    \end{tabular}
    \caption{Benchmark II: The computed displacements and the corresponding von Mises stress over the diaphragm. 
    This solution was obtained using the unfitted RBF-FD discretization with internodal distance $0.006$, an oversampling parameter 
    $q=5$ and an appended polynomial basis of degree $p=5$. \victor{Due to the scaling applied to the geometry (see Section~\ref{sec:smooth}), 
    the displayed results for displacement and stress should be multiplied with $s_{\Omega}^{-1}=0.15692$ m to recover the results corresponding 
    to the unscaled problem.}}
    \label{fig:experiments:Diaphragm:Robin:solution}
\end{figure}

The solution is given in Figure \ref{fig:experiments:Diaphragm:Robin:solution}, where we can see the spatial distribution of 
displacements and the Von Mises stress. To obtain this figure we used $h=0.006$, $p=5$, and $q=5$. We observe that 
there is a slight thickening in regions 2, 3, 7 and 8, according to the labels from Figure \ref{fig:smoothing:bcPosition}, which 
corresponds to the imposed thickening up to some extent. In regions 4 and 6, we observe a slight bend towards the interior, and in 
region 5, a translation in the negative vertical direction. This behavior does not entirely mimic the physiological contraction of 
the diaphragm, 
however, constructing a more accurate model is planned as future work. 
This solution serves as a test to understand whether the
unfitted RBF-FD method can handle the problem with the smoothed Robin boundary condition, 
as well as to explore how this type of boundary condition affects the behaviour of the solution. 
    
\subsection{Convergence under node refinement}
We validate the numerical solution by studying convergence under node refinement. 
Our choice of reference solution for this study is a fine solution computed using the unfitted RBF-FD method with an internodal distance 
$h=0.002$ and a polynomial degree $p=5$ used for constructing the local approximations. 
The relative error against the reference solution is computed using~\eqref{eq:dirichlet:error}. 
The relation between the considered $h$ and the number of points $N$ in $X$ is given in Table~\ref{fig:experiments:Dirichlet:Diaphragm:h_to_N_table}.

\begin{figure}[h!]
    \centering
    \begin{tabular}{ccc}
        \hspace{0.8cm}\textbf{Displacement }$\mathbf{\tilde u_1}$ & \hspace{0.8cm}\textbf{Displacement }$\mathbf{\tilde u_2}$ & \hspace{0.8cm}\textbf{Von Mises} \\
        \includegraphics[width=0.3\linewidth]{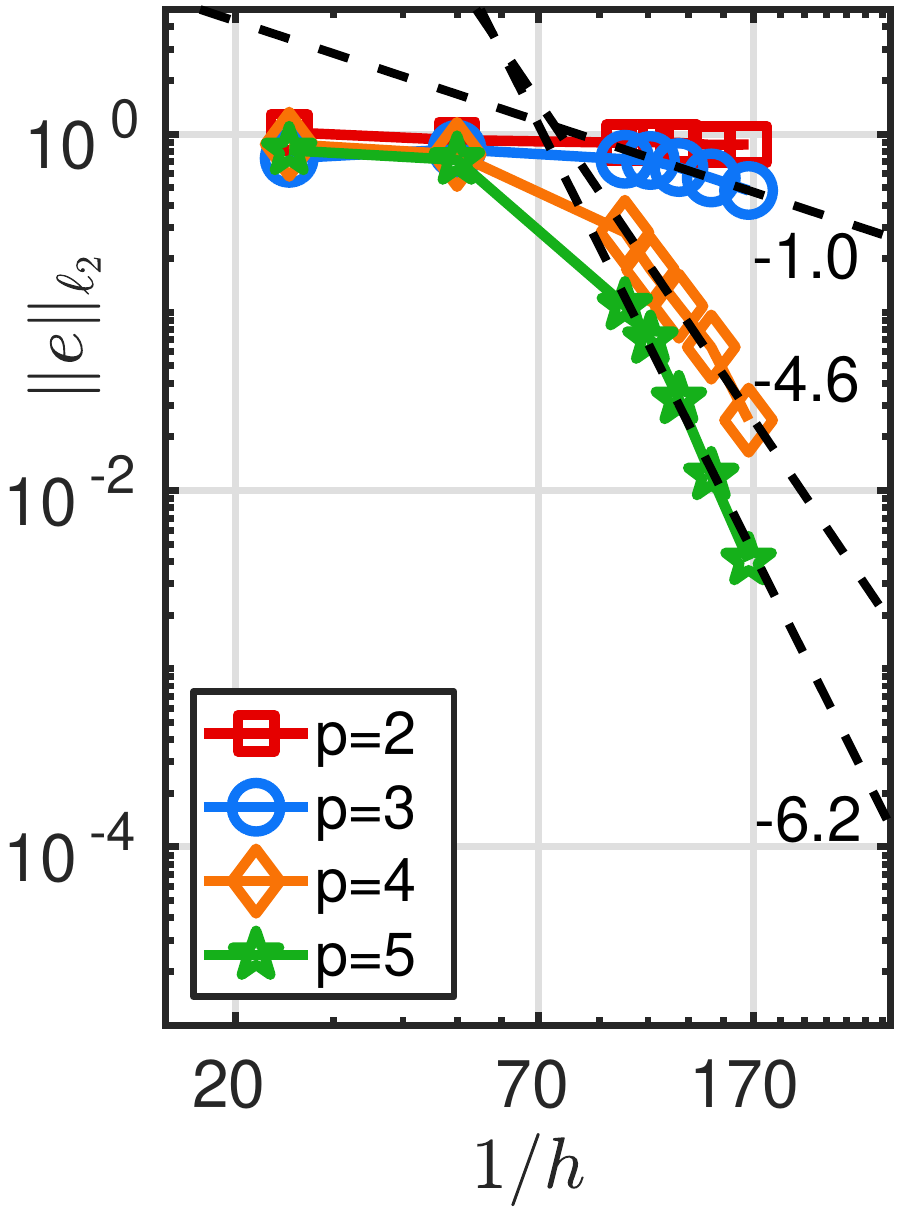} &  
        \includegraphics[width=0.3\linewidth]{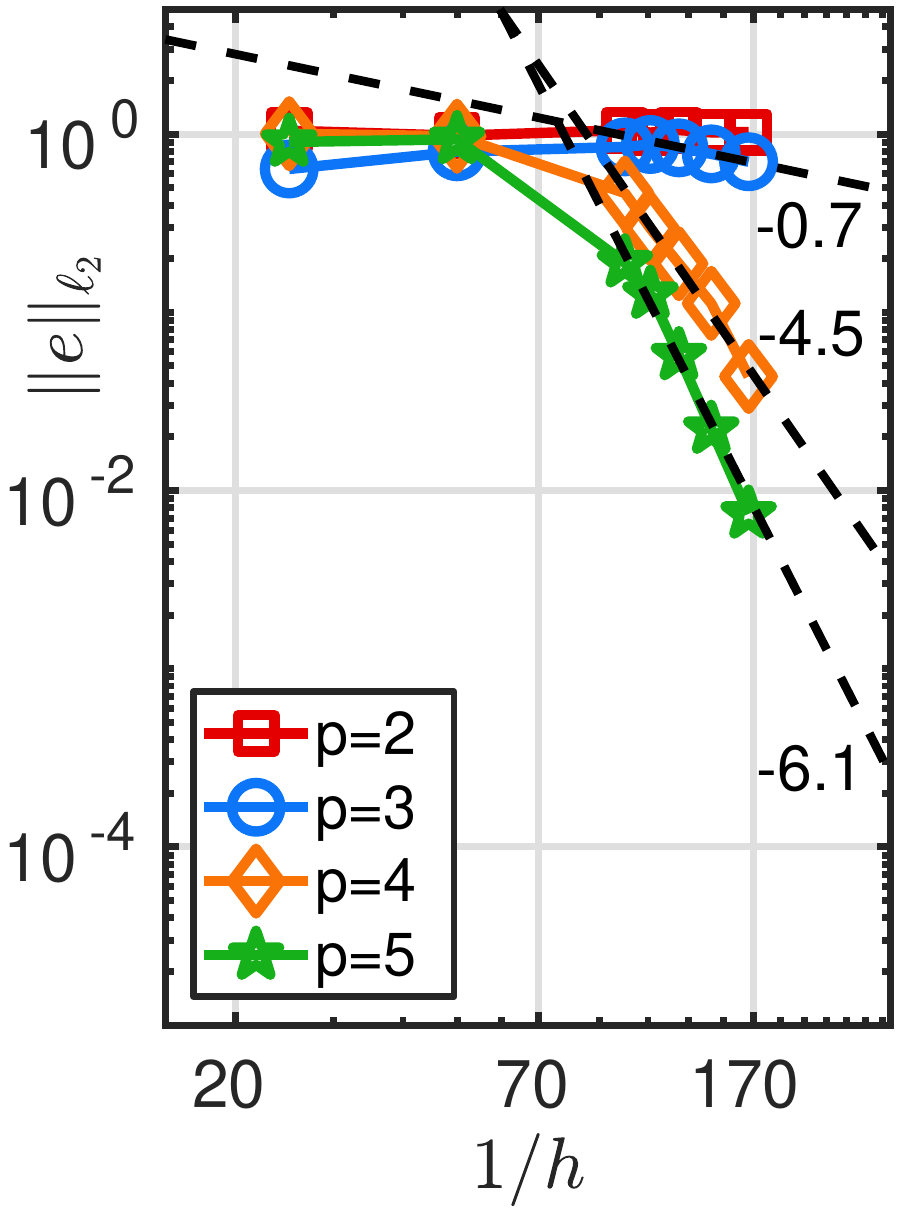} &  
        \includegraphics[width=0.3\linewidth]{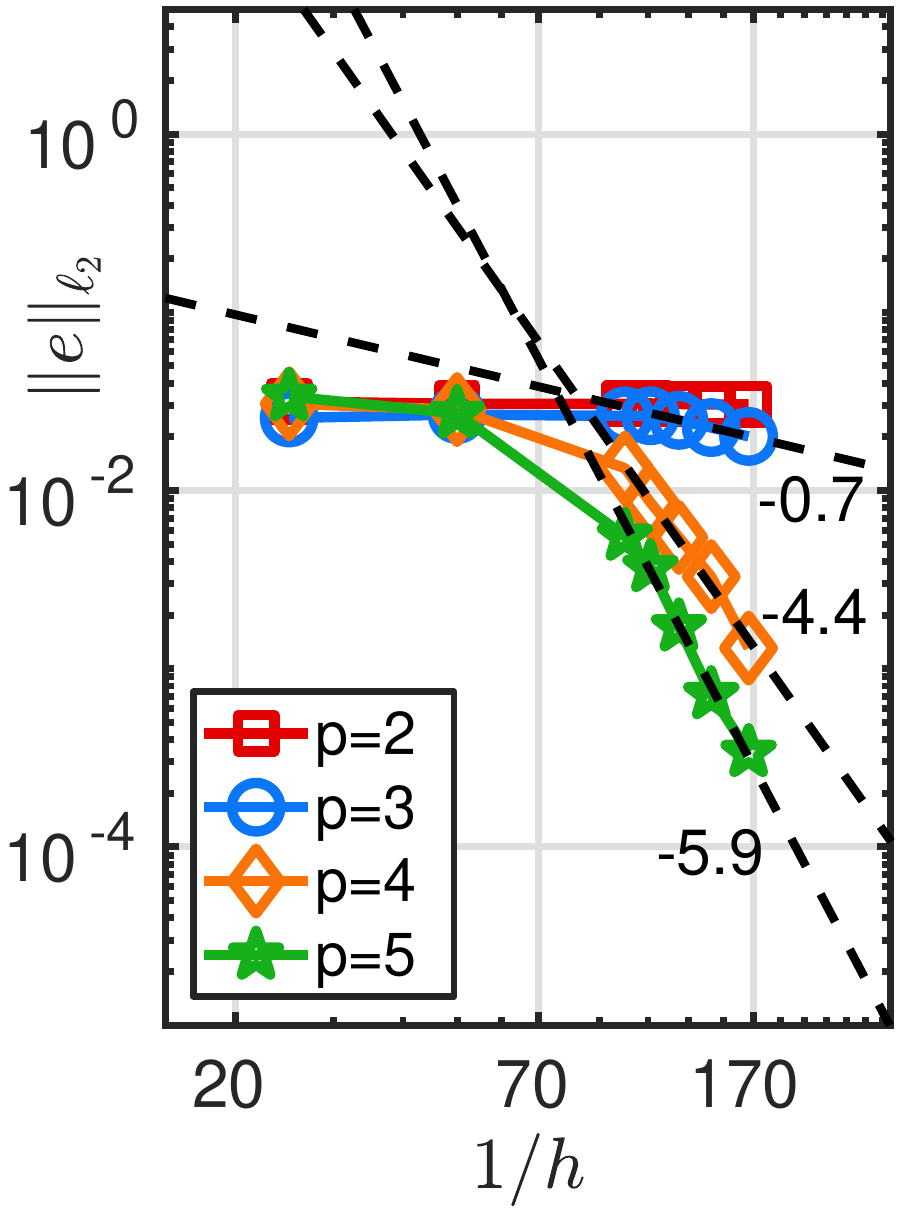}  
    \end{tabular}
    \caption{Benchmark II: Convergence of the displacements $\tilde u_1$ and $\tilde u_2$ and the Von Mises stress for different polynomial degrees $p$, 
    towards a dense numerical reference computed using the unfitted RBF-FD-LS method with $h=0.002$, $p=5$.}
    \label{fig:experiments:Diaphragm:Robin:convergence_selfRef}
\end{figure}

The errors for different choices of the polynomial degree $p$ are displayed in Figure \ref{fig:experiments:Diaphragm:Robin:convergence_selfRef}. 
\elisabeth{We do not observe convergence for any $p$ when $h$ is too large.} The reason for this is that the problem is not resolved yet. 
When $h$ is sufficiently small, and when $p=4$ or $p=5$, the solution converges at least with order $p-1$, which is desired. When $p=2$, 
we do not see convergence, since we would need an even higher resolution for this (small) polynomial degree. When $p=3$ we observe an approximately 
first order convergence for a sufficiently small $h$. This order is expected to increase to $p-1=2$ if $h$ is refined further. Low-order or no convergence 
when $p$ is small advocates using a higher order method.

\begin{figure}[h!]
    \centering
    
    \begin{tabular}{cc}
    \textbf{Displacement $\mathbf{\tilde u_1}$} & \textbf{Displacement $\mathbf{\tilde u_2}$} \\
    \includegraphics[width=0.49\linewidth]{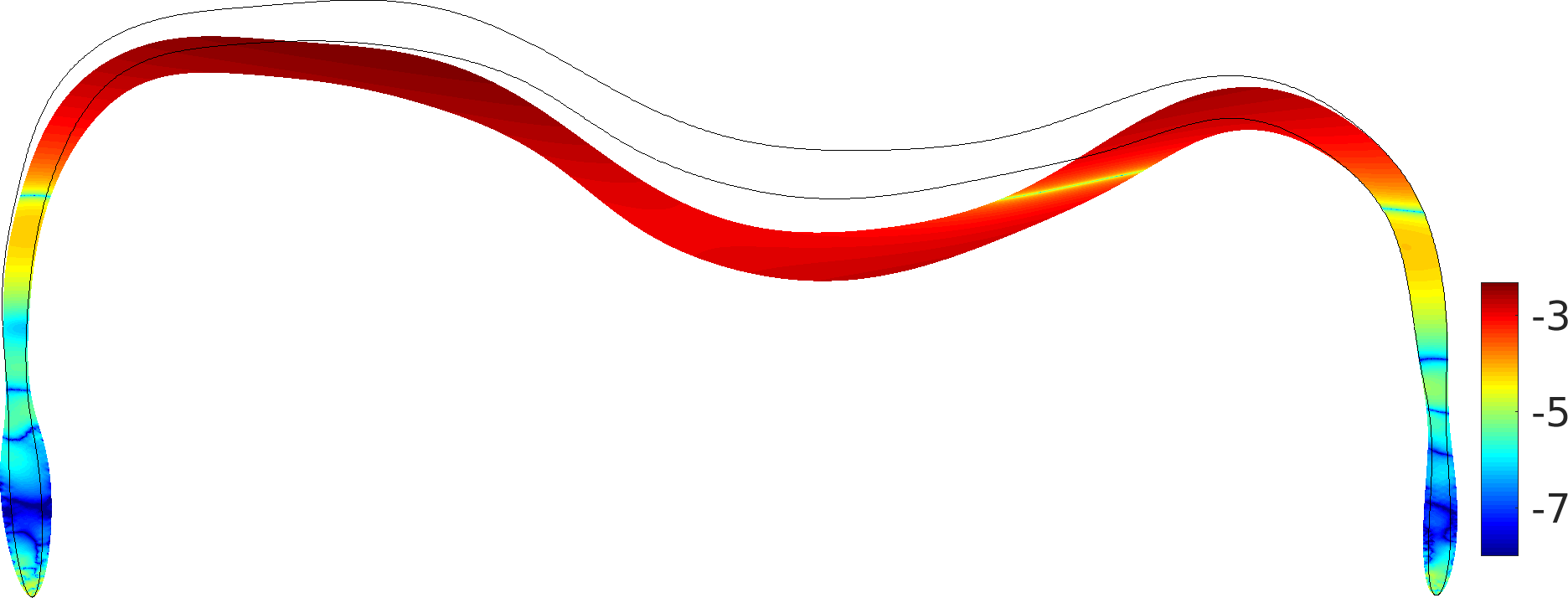} &
    \includegraphics[width=0.49\linewidth]{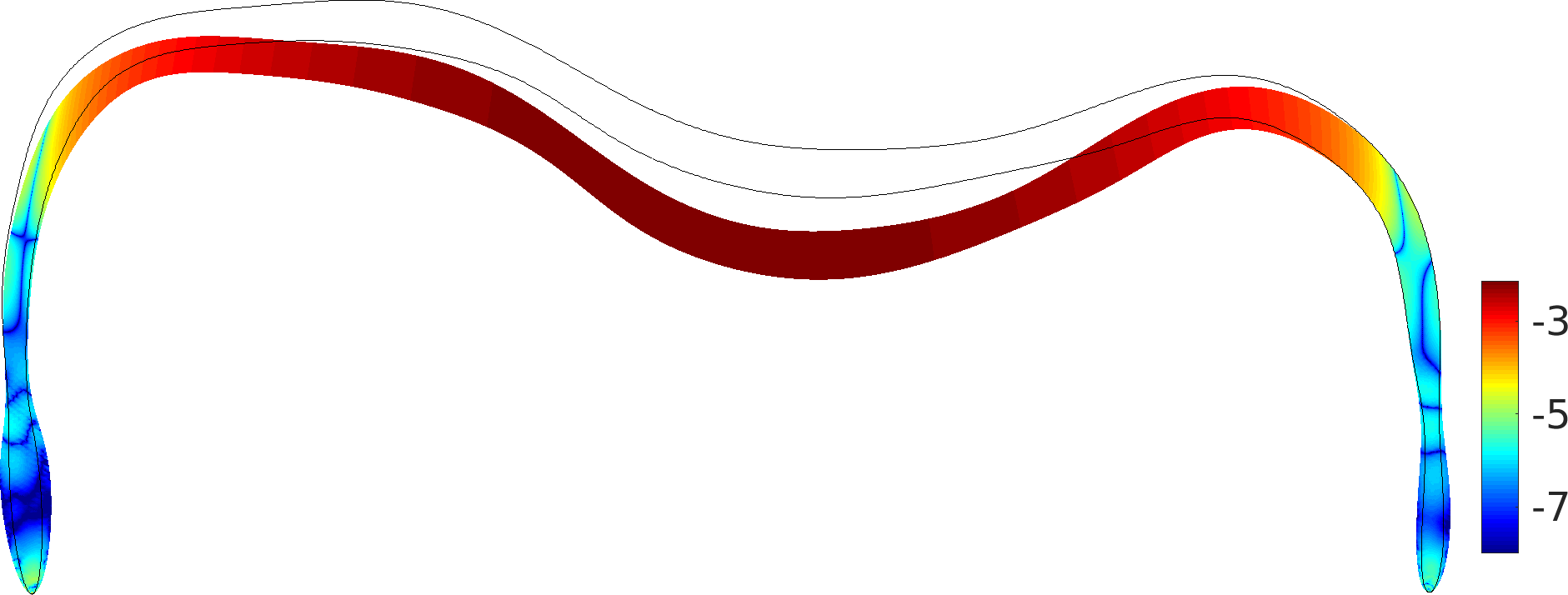} \\
    \multicolumn{2}{c}{\textbf{Von Mises stress}} \\
    \multicolumn{2}{c}{\includegraphics[width=0.6\linewidth]{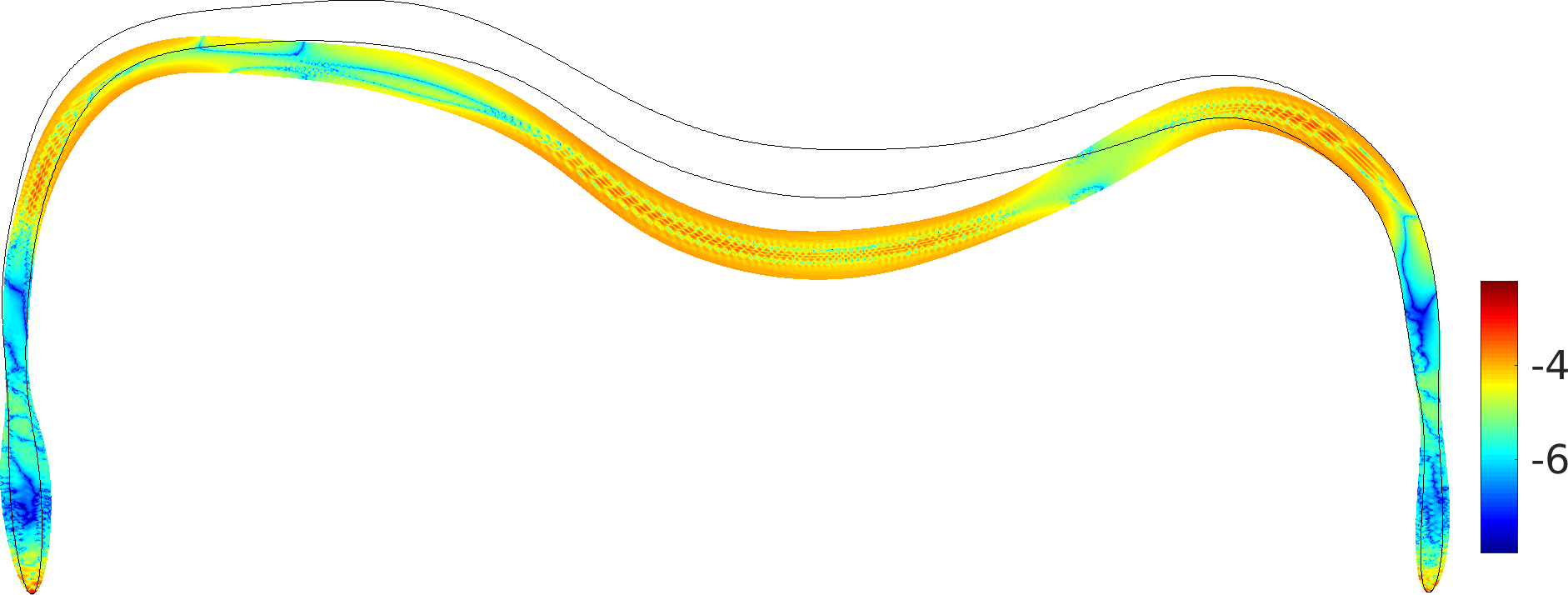}}
    \end{tabular}
\caption{Benchmark II: Error distribution in logarithmic scale for $h=0.006$ and $p=5$ when a fine unfitted RBF-FD-LS solution ($h=0.002$, $p=5$) 
         is taken as reference for computing the error.}
    \label{fig:experiments:Diaphragm:Robin:error_spatial_selfRef}
\end{figure}

The spatial distribution of the error computed using $h=0.006$ and $p=5$ is given in Figure
\ref{fig:experiments:Diaphragm:Robin:error_spatial_selfRef}. Here the error is largest close to 
region 5 from Figure \ref{fig:smoothing:bcPosition}, where we have enforced the traction boundary 
condition, which includes derivatives. In our experience the error is normally larger in the regions 
where derivative boundary conditions are imposed compared with the regions where a Dirichlet condition is 
imposed, and we also observed this in \cite{tominec2020unfitted,ToLaHe20}.

\section{Conclusions}
\label{sec:conclusions}

The unfitted RBF-FD method in the least-squares setting \cite{tominec2020unfitted} provided a robust framework 
for solving the linear elasticity system of PDEs over the simplified diaphragm geometry. 

We could also use the unfitted RBF-FD method to smooth the geometry curve and the boundary data. 
Ensuring that all components of the model were smooth allowed us to achieve high-order convergence, which reduces the number of unknowns 
needed for a given approximation error. In the experiments (not unexpectedly), Benchmark II, with Robin boundary conditions, proved more 
challenging to solve than Benchmark I with Dirichlet conditions. We needed higher resolution to achieve the same error level, and we did 
not see any convergence for the larger values of $h$. In the error plots, we could also see that the error is largest in the area where only 
the traction condition is active.
%
%
%
In the present work, we manufactured the data for the boundary conditions, but the aim is to eventually use measured pressure values 
and the solid body rotation of the ribs as boundary input data. This is most similar to the more challenging Benchmark II. 
There are a large number of transition zones, where the type of boundary condition changes, while we still expect a smooth 
behaviour of the solution. We were able to achieve a smooth solution and high-order convergence by using the proposed smoothing approach, 
and we will use that for the real application. Further investigations are needed regarding how to choose smoothing parameters such as the size of the transition zone.

To validate the unfitted RBF-FD solver, we measured convergence both against a self reference and against a highly resolved linear FEM solution. 
The results show that the two methods agree to high accuracy. 

Future work includes employing an unfitted RBF-FD-LS method to solve a more complex elastic problem, 
where we are going to use the 3D geometry and forcing data, both extracted from the CT (computed tomography) images of the diaphragm.

\section*{Acknowledgments}
The INVIVE project is funded by the Swedish Research Council, grant number 2016-04849.

\bibliographystyle{elsarticle-num}
\bibliography{refs}

\end{document}